\documentclass{amsart}
\usepackage{enumitem}
\usepackage[all]{xy}
\usepackage[margin=1.00in, marginparwidth=1.5cm, marginparsep=0.5cm]{geometry}

 \usepackage{booktabs}
 \usepackage{siunitx}

\usepackage[table]{xcolor}

\usepackage{lipsum} 

\usepackage{verbatim}

\usepackage{mathtools}

\usepackage{bbold}

\usepackage{placeins}

\usepackage{amssymb}
\usepackage{array}
\usepackage{xcolor}

\newcolumntype{C}[1]{>{\centering\arraybackslash}p{#1}}

\usepackage{mathrsfs}
\usepackage{physics}
\usepackage{amsthm}
\usepackage{amsmath,amsfonts,graphicx,mathdots,cite,wrapfig}
\usepackage[colorlinks,linkcolor=blue,citecolor=blue,pagebackref,hypertexnames=false, breaklinks]{hyperref}
\usepackage{cite}
\usepackage{mwe}
\usepackage{caption}
\usepackage{caption} 
\usepackage{subfigure} 
\usepackage{marginnote}

\newtheorem{theorem}{Theorem}[section]
\newtheorem{lemma}[theorem]{Lemma}
\newtheorem{proposition}[theorem]{Proposition}
\newtheorem{corollary}[theorem]{Corollary}

\theoremstyle{definition}
\newtheorem{definition}[theorem]{Definition}

\theoremstyle{remark}
\newtheorem{remark}[theorem]{Remark}

\numberwithin{equation}{section}

\DeclareMathOperator{\integers}{{\mathbb{Z}}}
\DeclareMathOperator{\nn}{\mathbb{N}}
\DeclareMathOperator{\rr}{\mathbb{R}}
\DeclareMathOperator{\complex}{\mathbb{C}}
\DeclareMathOperator{\hamilton}{\mathbb{{\textbf{H}}}}

\DeclareMathOperator{\jacobi}{\textbf{J}}

\newcommand{\bigdots}{\mbox{\normalfont\Large\bfseries $\ddots$}}
\newcommand{\DiriAMO}{\mbox{\normalfont\Large\bfseries $\hamilton_{p,\beta, \alpha, \theta }^{(l), D}$}}
\newcommand{\DiriAMOn}{\mbox{\normalfont\Large\bfseries $\hamilton_{p,\beta, \frac{k}{3^n} , 0 }^{(n), D}$}}

\begin{document}



\title{Spectral decimation of a self-similar version of almost Mathieu-type operators}


\author{Gamal Mograby}
\address{Mathematics Department, University of Maryland, College Park, MD 20742-4015 , USA}
\email{gmograby@umd.edu, gamal.mograby@uconn.edu}

\author{Radhakrishnan Balu}
\address{Radhakrishnan Balu, Department of Mathematics \& Norbert Wiener Center for
	Harmonic Analysis and Applications, University of Maryland, College Park, MD 20742, USA}
\email{radhakrishnan.balu.civ@mail.mil}

\author{Kasso A. Okoudjou}
\address{Kasso A. Okoudjou, Mathematics Department, Tufts University, Medford, MA 02155, USA}
\email{kasso.okoudjou@tufts.edu}

\author{Alexander Teplyaev}
\address{Alexander Teplyaev, Mathematics \& Physics Department, University of Connecticut, Storrs, CT 06269, USA}
\email{alexander.teplyaev@uconn.edu}


\date{\today}


\begin{abstract}
We introduce and study self-similar versions of the one-dimensional almost Mathieu operators. Our definition is based on a class of self-similar Laplacians instead of the standard discrete Laplacian, and  includes the classical almost Mathieu operators as a particular case. Our main result establishes that the spectra of these  self-similar almost Mathieu operators can be described by the spectra of the corresponding  self-similar Laplacians through the spectral decimation framework used in the context of spectral analysis on fractals. The spectral type of the self-similar Laplacians used in our model are singularly continuous. The self-similar almost Mathieu operators also have singularly continuous spectrum for specific parameters. In addition, we derive an explicit formula of the integrated density of states of the self-similar almost Mathieu operators as the weighted pre-images of the balanced invariant measure on a specific Julia set.
\end{abstract}

\maketitle

\tableofcontents

 \section{Introduction}\label{sec1}
 
 The investigation of the properties of quasi-periodic Schr\"odinger-type operators remains very active drawing techniques from different areas of mathematics and physics \cite{marx_jitomirskaya_2017, Wilkinson2017,Akkermans2021}.  
The special case of the almost Mathieu operators (AMO)  can be traced back to Harper who proposed a model to describe crystal electrons in a uniform magnetic field \cite{Harper_1955}. Subsequently, Hofstadter observed that the spectra of the AMO can be fractal sets \cite{Hofstadter1976}. We refer to \cite{dinaburg_one-dimensional_1976, moser_example_1981} for more early examples of such operators whose spectra are Cantor-like sets, and to 
\cite{avila_ten_2009, bellissard_cantor_1982, jitomirskaya_metal-insulator_1999,  v_mouche_coexistence_1989} for more results on the AMO and references therein. 

Independently,  a line of investigations of self-similar Laplacian operators on graphs, fractals, and networks has emerged \cite{Rammal1984SpectrumOH, RammalToulouse1983, Alexander1984,KadanoffAlexander1983}. A fundamental tool in this framework is the spectral decimation method, initially  used in physics to compute the spectrum of the Laplacian on a Sierpinski lattice \cite{KigamiAnaOnFractalsBook, BellissardRenormalizationGroup1992,StrichartzTep2012,Strichartz2003FractafoldsBO,malozemovteplyaev2003}. At the heart of this method is the fact that  the spectrum of this Laplacian  is completely described in terms of iterations of a rational function. For an overview of the modern mathematical approaches, applications, and extensions of the spectral decimation methods we refer to  \cite{Shirai2000, StrichartzTrafoGraphLap2010 ,StrichartzMethodofAver2001, Fukushima1992OnAS, ShimaSierpinski1993,  ShimapreSierpinski1991, TeplyaevInfiniteSG1998,  BobsBook, BajorinSteinhurst2008, kron_asymptotics_2003} and references therein.

Recently, Chen and Teplyaev \cite{ChenTeplyPQmodel2016} used the general framework of the spectral decimation method to investigate the appearnce of the singular continuous spectrum of a family of Laplacian operators. One of the   ideas used in establishing this result is that these Laplacians are naturally related to self-similar operators  with corresponding self-similar structures \cite{malozemovteplyaev2003} which allows to use complex dynamics techniques. 
\begin{figure} 
  \begin{minipage}[b]{0.5\linewidth}
    \hspace*{-1cm}
    \centering
    \includegraphics[width=1.2\linewidth]{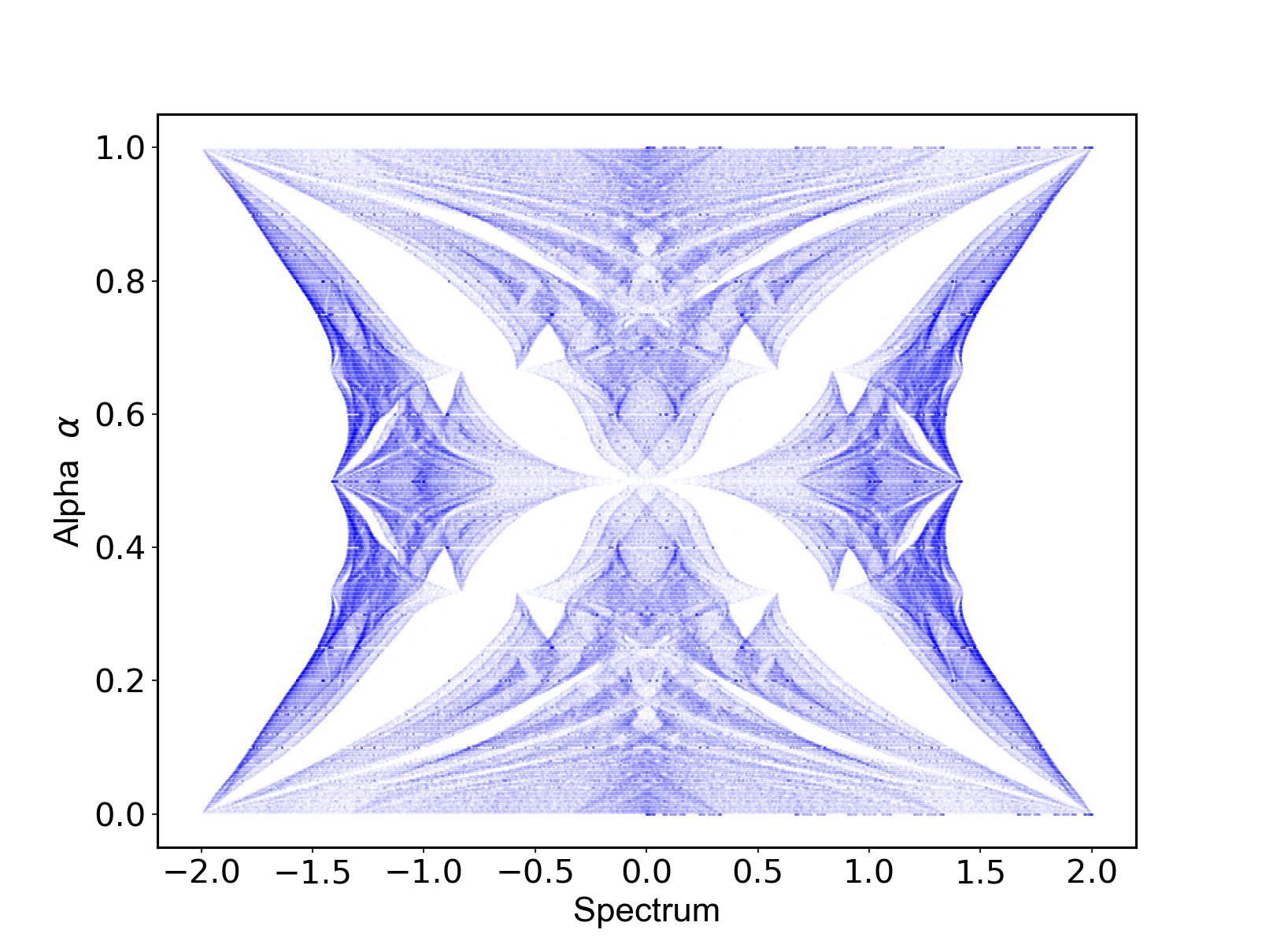} 
  \end{minipage}
  \begin{minipage}[b]{0.5\linewidth}
    \centering
    \includegraphics[width=1.16\linewidth]{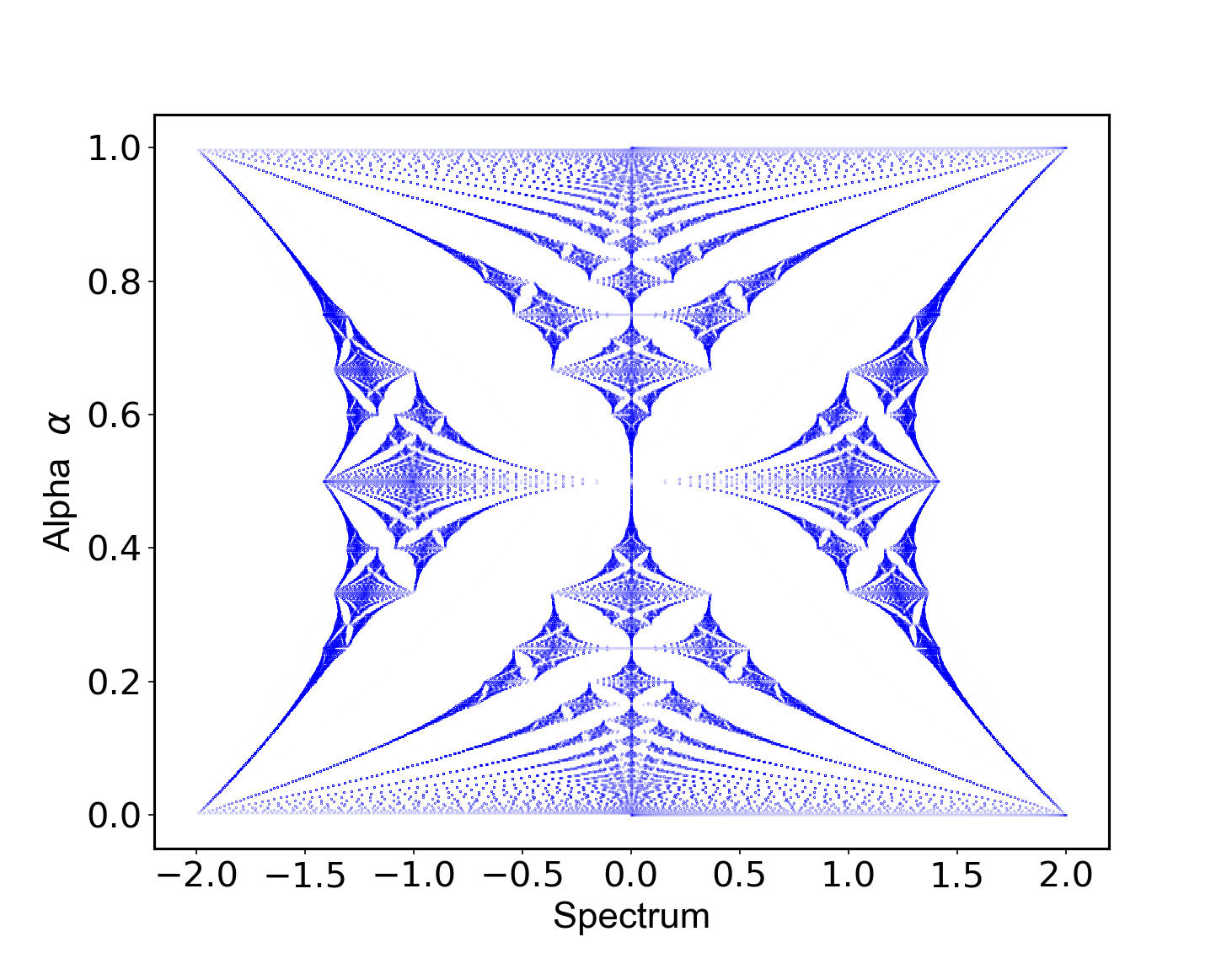} 
  \end{minipage} 
  \caption{The left panel of the figure shows a Hofstadter butterfly for a self-similar almost Mathieu operator corresponding to $\frac{1}{3}$-Laplacian whose spectrum is a Cantor set.  For comparison, the (classical) Hofstadter butterfly corresponding to  the standard  AMO is shown in right panel.}
  \label{fig:HofstadterButterfliesVaryDiffPara}
\end{figure}

\begin{figure} 
  \begin{minipage}[b]{0.5\linewidth}
    \hspace*{-1cm}
    \centering
    \includegraphics[width=1.16\linewidth]{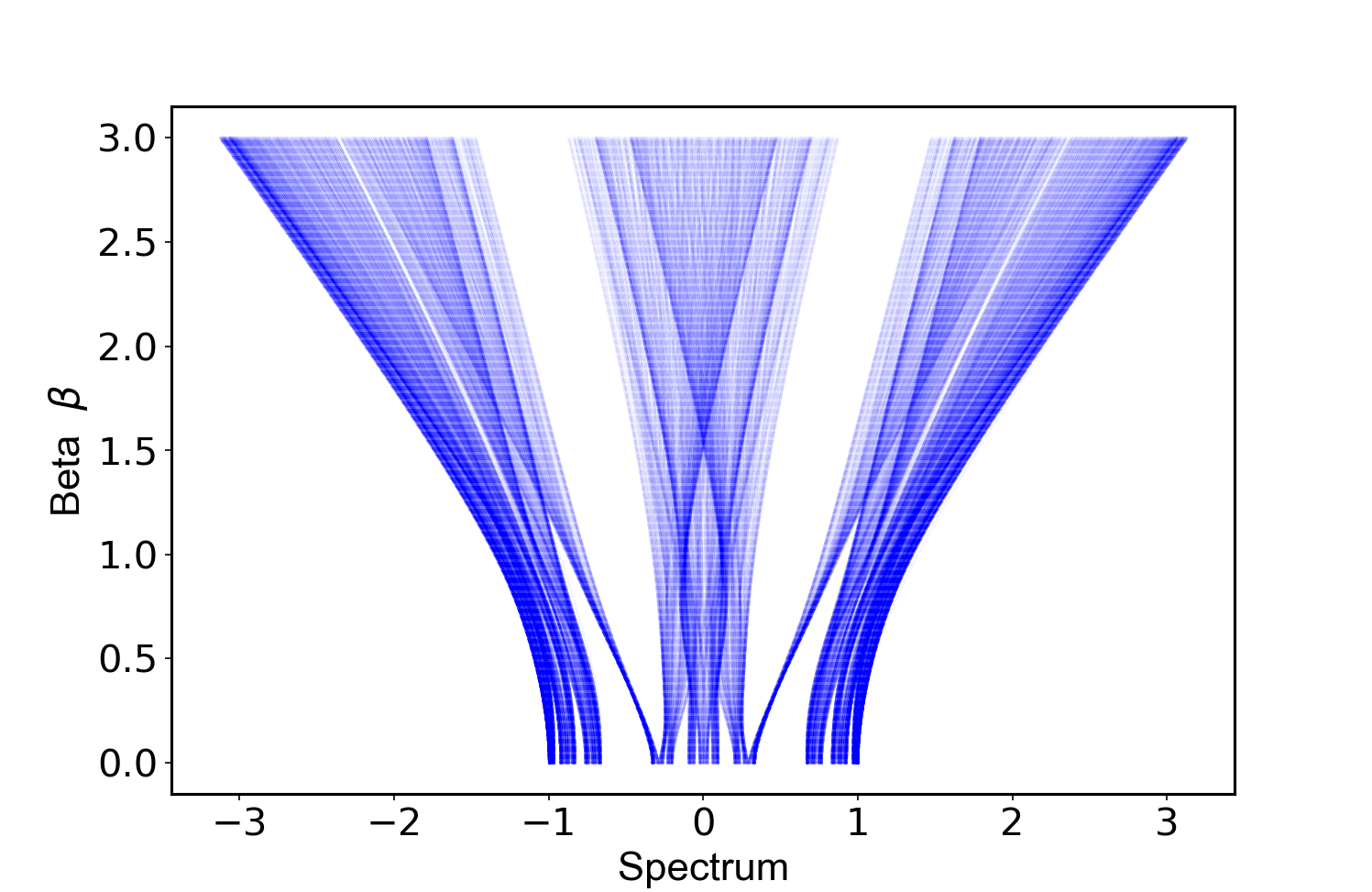} 
  \end{minipage}
  \begin{minipage}[b]{0.5\linewidth}
    \centering
    \includegraphics[width=1.2\linewidth]{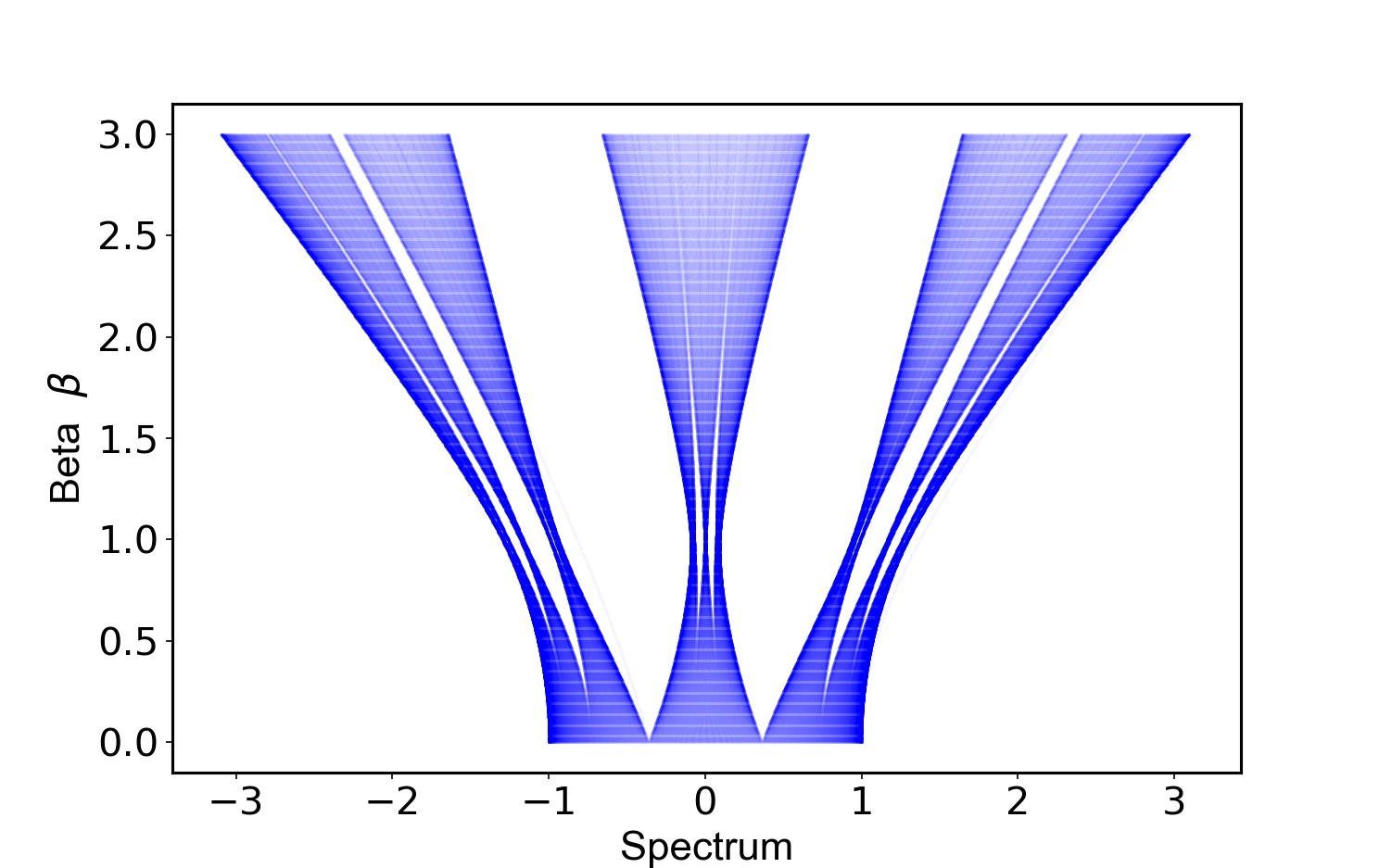} 
  \end{minipage} 
  \caption{ The spectrum of $\hamilton_{p,\beta, \alpha, \theta}$ is plotted ($x$-axis) for the fixed parameters $\alpha = \frac{\sqrt{5}-1}{2}$ and $\theta=0$ while varying $\beta \in [0,3]$ (y-axis). The parameter $p$ is equal to $\frac{1}{3}$ for the left- and equal to $\frac{1}{2}$ for the right-panel. Both panels depict examples corresponding to an irrational $\alpha$. }
  \label{fig:VaryBetaIrrationalAlpha}
\end{figure}
The present paper is a first in what we expect to be a research program dealing with quasi-periodic Schr\"odinger-type operators  on self-similar sets such as fractals and graphs.  Our goal is to initiate the study of  a generalization of the discrete almost Mathieu operators for self-similar situations. 
In this paper we begin by considering  finite or half-integer one-dimensional lattices endowed with   particular self-similar structures. 
More general Jacobi matrices  will be considered in the forthcoming work  \cite{BaluMograbyOkoudjouTeplyaevJacobi2021}.

In this setting, the \textit{self-similar almost Mathieu operators (s-AMO)} are formally defined in~\eqref{eq:AMOversionPQ}  and will be denoted  by $\hamilton_{p,\beta, \alpha, \theta }$, for  $\alpha \in \rr$, $\theta \in [0,2 \pi)$ and $\beta \in \rr$.
 As we will show, these operators can be viewed  as limits of finite dimensional analogues that can be completely understood using the spectral decimation methods developed by Malozemov and Teplyaev \cite{malozemovteplyaev2003}. Furthermore, the s-AMO we consider, are defined in terms of self-similar Laplacians $\{\Delta_p\}_{p \in (0,1)}$ which are given by~\eqref{pLaplace}. This class of self-similar Laplacians was first investigated in \cite{TeplyaevSpectralZeta2007} and arises naturally when studying the unit-interval endowed with a particular fractal measure, see also the related work \cite{fang_spectral_2019, Derfel_2012,chan_one-dimensional_2015,Bird_Ngai_Tep_2003}.
 Moreover, when $p=\frac{1}{2}$, then the self-similar almost Mathieu operators coincide up to a multiplicative constant with the standard one-dimensional almost Mathieu operators (see (\ref{eq:UsualAMO})). Whereas, in the case of AMO when the magnetic flux is a fraction the resulting spectrum is absolute continuous,  the fractal case 
 displays a spectrum that is singular continuous.

 The paper is organized as follows. In Section~\ref{HierarchicalAMO} we introduce the notations and the definition of the self-similar structure we impose on the half-integer lattice.  In the first part of Section~\ref{sec-spectralDeci}, we focus on the discrete and finite s-AMOs and describe their spectra using the spectral decimation method, see Section~\ref{sec-finitegraphs}. Subsequently,  in Section~\ref{sec:infinitegraphs} we prove one of our main results, Theorem~\ref{thm:SpectralSimiAMOandpqModelnew}, which states that the spectra of the AMO on $\integers_+$ can be completely described using the spectral decimation method when the parameter $\alpha$ belongs to a dense set of numbers.  Moreover,  these operators have purely singularly continuous spectra when $p\neq 1/2$.  As will be seen, Theorem \ref{thm:SpectralSimiAMOandpqModelnew}  provides a useful algebraic tool to relate the spectra of the almost Mathieu operators  $\hamilton_{p,\beta, \alpha, \theta }$ to that of the family of self-similar Laplacians $\Delta_p$. In Section~\ref{IDSsection},  using the fact that the spectrum of a self-similar Laplacian $\Delta_{p}$ is the Julia set $\mathcal{J}(R_{\Delta_{p}})$ of a polynomial (defined in ~\eqref{eq:specDeciFuncLap}), we derive in Theorem \ref{thm:IDSofAMO} an explicit formula for the density of states of  $\hamilton_{p,\beta, \frac{k}{3^n}, 0 }$ by identifying it with the weighted pre-images of the balanced invariant measure on the Julia set  $\mathcal{J}(R_{\Delta_{p}})$.  As a corollary, we obtain a gaps labeling statement for  $\hamilton_{p,\beta, \frac{k}{3^n}, 0 }$. In Section~\ref{sec:examplesappl} we present some numerical simulations pertaining the spectra of the s-AMO as well as the integrated density of states for a variety of parameters.

A first illustration of our numerical results is Figure~\ref{fig:HofstadterButterfliesVaryDiffPara}. The left panel of the figure shows a Hofstadter butterfly for a self-similar almost Mathieu operator corresponding to $\frac{1}{3}$-Laplacian.  For comparison, the Hofstadter butterfly corresponding to  the standard  AMO is shown in the right panel. In both cases, the spectrum is plotted ($x$-axis) for the fixed parameters  $\beta=1$ and $\theta=0$ while varying $\alpha \in [0,1]$ (y-axis). For all $\alpha \in \{\tfrac{k}{3^n}, \, k=0, 1, \hdots, 3^n-1\}_{n=1}^l$ where $l\geq 1$, Theorem \ref{thm:SpectralSimiAMOandpqModelnew} describes the difference in these two figures as a transformation given by a spectral decimation function. Note that in the  standard  AMO case, corresponding to $p=1/2$ in our formulation, many important results are also obtained  for $\alpha$ irrational, but in the fractal setting the methods for irrational $\alpha$ are not developed yet.  Figure~\ref{fig:VaryBetaIrrationalAlpha} depict examples corresponding to the case of an irrational $\alpha = \frac{\sqrt{5}-1}{2}$.  The spectrum of $\hamilton_{p,\beta, \alpha, \theta }^{(l)}$ is plotted ($x$-axis) for the fixed parameters $\alpha = \frac{\sqrt{5}-1}{2}$ and $\theta=0$ while varying $\beta \in [0,3]$ (y-axis).

We end this introduction with a perspective on a general framework underlying the present paper. 
 In a  forthcoming and companion paper \cite{BaluMograbyOkoudjouTeplyaevJacobi2021}, we  identify a class of Jacobi operators that extend the present results to almost Mathieu operators  defined in the fractal setting. We refer to this class of operators as piecewise centrosymmetric Jacobi operators  \cite{centrosymmetric1,centrosymmetric2,centrosymmetric3}. 
 In this general setting we show that the spectral decimation function arises from a  particular system of orthogonal polynomials associated with the aforementioned Jacobi matrix. In particular, this spectral decimation function is computable using a three-term recursion formula associated to this system of orthogonal polynomials. In the process, we avoid the Schur complement computation that could involve resolvent calculations of large matrices. Additionally, the general setting we consider in \cite{BaluMograbyOkoudjouTeplyaevJacobi2021}
can be further extended to higher-dimensional graphs and would allow us to define Jacobi-type operators on graphs like the Sierpinski lattices or Diamond graphs. We plan to use this approach to investigate some of the questions in \cite{AvniSimon2020}.

\section{Self-similar Laplacians and almost Mathieu operators} \label{HierarchicalAMO}
In this section we introduce the notations and the definition of the self-similar structure we impose on the half-integer lattice. This self-similar structure describes a random walk on the half-line and gives rise to a class of self-similar probabilistic graph Laplacians $\Delta_p$. Moreover, it provides a natural finite graph approximation for the half-integer lattice. Regarding an almost Mathieu operator as a Schr\"odinger-type operator of the form $\Delta + U$ (where $U$ is a potential operator), allows us to define the class of self-similar almost Mathieu operators as $\Delta_p + U$.

\subsection{Self-similar $p$ Laplacians on the half-integer lattice}\label{subsec2.1}

 We consider a family of  self-similar Laplacians on the integers half-line. This class of Laplacians was first investigated in \cite{TeplyaevSpectralZeta2007} and arises naturally when studying the unit-interval endowed with a particular fractal measure. For more on this Laplacian and some related work we refer to  \cite{Derfel_2012, chan_one-dimensional_2015, wave17}. The Laplacian's spectral-type was investigated in \cite{ChenTeplyPQmodel2016},  where the emerging of singularly continuous spectra was proved. Furthermore, this class of Laplacians serves as a toy model for generating singularly continuous spectra. In this section, we   introduce the $p$-Laplacians and review some of its properties that will be needed to state and prove our results, and refer to  \cite{ChenTeplyPQmodel2016} for more details. We also introduce a corresponding self-similar structure on the half-integer line.

Let $\mathbb{Z}_+$ be the set of nonnegative integers and $\ell(\mathbb{Z}_+)$ be the linear space of complex-valued sequences $(f(x))_{x \in \mathbb{Z}_+}$. Let $p\in (0,1)$, for each $x\in \mathbb{Z}_+ \setminus \{0\}$, we define $m(x)$ to be the largest natural number $m$ such that $3^m$ divides $x$. For $f \in \ell(\mathbb{Z}_+)$ we define a \textit{self-similar Laplacian} $\Delta_p$ by,
\begin{align}\label{pLaplace}
	(\Delta_p f)(x) = \left\{\begin{array}{ll} f(0)-f(1), & \text{if}~x=0 \\
		f(x)-(1-p)f(x-1)-pf(x+1), &\text{if}~3^{-m(x)}x \equiv 1~\pmod 3 \\
		f(x) - pf(x-1)-(1-p)f(x+1), &\text{if}~3^{-m(x)}x \equiv 2~\pmod 3
	\end{array}\right..
\end{align}
We equip $\ell(\mathbb{Z}_+)$  with its canonical basis  $\{\delta_x\}_{x \in \mathbb{Z}_{+}}$ where 
\begin{equation}
\label{eq:canonicalBasis}
\delta_x(y)
 =
  \begin{cases}
  \  0     & \quad \text{if } x \neq y \\
  \  1 & \quad \text{if } x=y.
  \end{cases}
\end{equation}
 The matrix representation of $\Delta_p$  with respect to the canonical basis has the following Jacobi matrix

\begin{equation}
\jacobi_{+,p} =
\begin{pmatrix}
1 & -1 & 0 & 0 & 0 & 0 & 0 & 0 & \dots  \\
p-1 & 1 & - p & 0 & 0 & 0 & 0 & 0 & \dots \\
0 & - p & 1 & p-1 & 0 & 0 & 0 & 0 & \dots \\
0 & 0 & p-1 & 1 & - p & 0 & 0 & 0 & \dots \\
0 & 0 & 0 & p-1 & 1 & - p & 0 & 0 & \dots \\
0 & 0 & 0 & 0 & - p & 1 & p-1 & 0 & \dots \\
0 & 0 & 0 & 0 & 0 & - p & 1 & p-1 & \dots \\
0 & 0 & 0 & 0 & 0 & 0 & p-1 & 1 &    \dots\\
\vdots & \vdots & \vdots & \vdots & \vdots & \vdots & \vdots & \vdots & \ddots
\end{pmatrix}.
\end{equation}
\text{ } \\
The case $p = \frac{1}{2}$ recovers the classical one-dimensional Laplacian (probabilistic graph Laplacian).  

We adopt the notation used to describe a random walk on the half-line with reflection at the
origin and refer to the off-diagonal entries in $\jacobi_{+,p}$ by the transition probabilities
\begin{equation}
\label{eq:transitionProbabilities}
p(x,y) = - \jacobi_{+,p}[x,y], \quad \text{ for } x \neq y.
\end{equation}
Let $\pi$ be a $\sigma$-finite measure on $\mathbb{Z}_+$. We define the Hilbert space
\begin{equation*}
    \ell^2(\integers_{+}, d \pi) = \{ \psi: \mathbb{Z}_+ \to \complex \ | \ \sum^{\infty}_{x =0} \ |\psi (x)|^2  \pi(x) < \infty  \}, \quad  \langle f, g \rangle_{\ell^2} = \sum_{x=0}^{\infty} \overline{f(x)} g(x) \pi(x).
\end{equation*}
Let  $n \in \mathbb{Z}_+$, the ($n$-th) Wronskian of $f,g \in \ell(\mathbb{Z}_+)$ is given by
\begin{align}
W_n(f,g)  =  \pi(n)  p(n,n+1)    \Big( \overline{f(n)}g (n+1)     -\overline{f(n+1)}g(n) \Big).
\end{align}
\begin{lemma}
\label{lem:SelfAdjointnessOfLap}
Let $f,g \in  \ell^2(\integers_{+}, d \pi)$ and $n \in \mathbb{Z}_+$. Assume that the measure $\pi$ satisfies the reversibility condition, i.e., $\pi(x)p(x,y)=\pi(y)p(y,x)$ holds for every  $x,y \in \mathbb{Z}_+$. Then the discrete Green's second identity holds. That is, we have: 
\begin{equation}
 \sum_{x=0}^n \overline{f(x)} \Delta_p g (x) \pi(x) -  \sum_{x=0}^n \overline{ \Delta_p f(x) } g (x) \pi(x) = W_n(f,g).
\end{equation}
Moreover,  the operator $\Delta_p$ is a bounded self-adjoint operator on $ \ell^2(\integers_{+}, d \pi)$.
\end{lemma}

\begin{proof}
Direct computation gives for $n \in \mathbb{Z}_+\backslash \{0\}$,
\begin{align*}
f(n)  \Delta_p g(n)   \pi(n) -   \Delta_p f (n)   g (n)  \pi(n) &=   W_n(f,g) - \pi(n)   p(n,n-1) \Big(f(n-1)g (n)    - f(n)   g(n-1)  \Big).
\end{align*}
Using the reversibility condition, i.e. $\pi(n) p(n,n-1) = \pi(n-1) p(n-1,n)$, we obtain 
\begin{align*}
f(n) \Delta_p g(n)   \pi(n) -  \Delta_p f (n)  g (n)  \pi(n) &=   W_n(f,g) - W_{n-1}(f,g).
\end{align*}
For $n=0$, we compute
\begin{align*}
f(0) \Delta_p g(0)   \pi(0) -  \Delta_p f (0)  g (0)  \pi(0) =   f(0) p(0,1) g(1) \pi(0)  - g(0) p(0,1) f(1) \pi(0) = W_0(f,g)
\end{align*}
Hence, a telescoping trick gives
\begin{equation*}
 \sum_{x=0}^n f(x) \Delta_p g (x) \pi(x) -  \sum_{x=0}^n  \Delta_p f(x)  g (x) \pi(x) = W_n(f,g).
 \end{equation*}
 For $f,g \in \ell^2( \mathbb{Z}_+, d \pi)$, we imply $ \langle f, \Delta_p g \rangle_{\ell^2} -  \langle \Delta_p f, g \rangle_{\ell^2} = \lim_{n \to \infty} W_n(f,g) =0$.
\end{proof}
\begin{figure}[!htb]
    \centering
  \includegraphics[width=1.\textwidth]{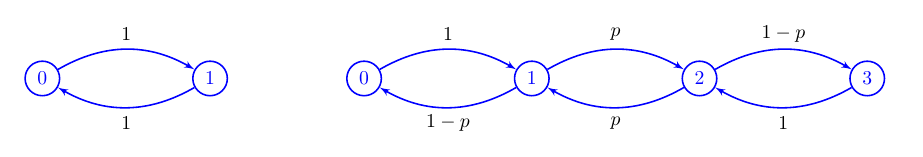}
    	\caption{(Left) Initializing the graph $G_0$. (Right) The graph $G_1$. While the vertices are labeled by the addresses, the labeling of the edges represents the transition probabilities \ref{eq:transitionProbabilities}. }
	\label{fig:G1_initialNeumann}
\end{figure}
\begin{figure}[htp]
\centering
 \includegraphics[width=1.\textwidth]{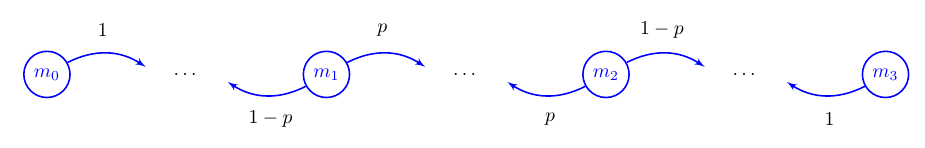}
\caption{The visual representation of the protograph indicates how to apply the substitution rule, see Definition \ref{def:finiteGraphApproxNeumann}.}
\label{fig:modelGraphNeumann}
\end{figure}
We regard the integers half-line $ \mathbb{Z}_+$ endowed with $\Delta_p$ as a hierarchical or substitution infinite graph, see \cite{malozemovteplyaev1995, malozemovteplyaev2003} for more details. We define a sequence of finite directed weighted graphs $\{G_l \}_{l \in \nn}$, such that $G_l = (V_l, E_l) $ is constructed inductively according to a substitution rule. We set $V_l =\mathbb{Z}_+ \cap [0,3^l]$ for all $l \geq 0$, where  $G_0 = (V_0, E_0)$ is the graph shown in Figure \ref{fig:G1_initialNeumann} (Left). We illustrate the substitution rule by constructing $G_1$ shown in Figure \ref{fig:G1_initialNeumann} (Right).  We first introduce the \textit{protograph} shown in Figure \ref{fig:modelGraphNeumann}, which consists of the four vertices $\{m_0, m_1,m_2, m_3\}$. We insert three copies of $G_0$ in the protograph according to the following rule. Between any two vertices $m_i$ and $m_{i+1}$, we substitute the three dots with a copy of $G_0$, identifying the vertex $0$ in $G_0$ with the vertex $m_i$, and the vertex $1$ in $G_0$ with the vertex $m_{i+1}$.  We substitute the edges $(0,1)$ and $(1,0)$ in $G_0$ with the corresponding directed weighted edges as indicated in the protograph, see Figure \ref{fig:modelGraphNeumann}. We repeat the procedure and insert copies of $G_0$ between the vertices, $m_0$, $m_1$, then  $m_1$, $m_2$ and finally $m_2$, $m_3$. The resulting linear directed weighted graph is denoted by $G_1$, Figure \ref{fig:G1_initialNeumann} (Right). The graph $G_1$ consists of $4$ vertices, which we rename to $\{0,1,2,3 \}$, so that $m_0$ corresponds to the vertex $0$, $m_1$ to $1$ , $m_2$ to $2$ and $m_3$ to $3$. In particular, this gives $V_1 =\mathbb{Z}_+ \cap [0,3^1]$ and $G_1$ can be viewed as a truncation of $\mathbb{Z}_+$ (regarded as a hierarchical infinite graph) to the vertices $\{0,1,2 ,3 \}$, whereby a reflecting boundary condition is imposed on the vertex $3$.
Similarly, we construct $G_2$ by inserting $G_1$ in the protograph, see Figure \ref{fig:G2Neumann}. 
\begin{definition}
\label{def:finiteGraphApproxNeumann}
 Let  $G_0 = (V_0, E_0)$ be the graph shown in Figure \ref{fig:G1_initialNeumann} (Left). We define the sequence of graphs $\{G_l \}_{l \in \nn}$ inductively. Suppose $G_{l-1}=(V_{l-1}, E_{l-1})$ is given for some integer $l \geq 1$, where $V_{l-1} =\mathbb{Z}_+ \cap [0,3^{l-1}]$. The graph $G_{l}=(V_{l}, E_{l})$ is constructed according to the following \textit{substitution rule}. We repeat the following steps for $i \in \{0,1,2\}$:
\begin{enumerate}
	\item  Insert a copy of $G_{l-1}$ between the two vertices $m_i$ and $m_{i+1}$ of the protograph shown in Figure \ref{fig:modelGraphNeumann} in the following sense. We identify the vertex $0$ in $G_{l-1}$ with the vertex $m_i$ and similarly, we identify the vertex $3^{l-1}$ in $G_{l-1}$ with the vertex $m_{i+1}$.
	\item  We substitute the edges $(0,1)$ and $(3^{l-1},3^{l-1}-1)$ in  $G_{l-1}$ with the corresponding directed weighted edges as indicated in the protograph, see Figure \ref{fig:modelGraphNeumann}.
\end{enumerate}
\end{definition}
The resulting linear directed weighted graph is denoted by $G_{l}=(V_{l}, E_{l})$. The graph $G_{l}$ consists of $3^l +1$ vertices, which we rename to $\{0,1,\dots,3^l \}$, so that $m_0$ corresponds to the vertex $0$, ... , $m_l$ corresponds to the vertex $3^l$. In particular, this gives $V_l =\mathbb{Z}_+ \cap [0,3^l]$. The vertices $0$ and $3^l$ are the boundary vertices of $G_l$, and we refer to them by $\partial G_l = \{0,3^l\}$. The interior vertices of $G_l$ are given by $V_l \backslash \partial G_l$. 

Each graph $G_{l}=(V_{l}, E_{l})$ is naturally associated with a \textit{probabilistic graph Laplacian}, denoted  $\Delta^{(l)}_{p}$, and given by 
\begin{equation*}
\Delta^{(l)}_{p} f (x)= \Delta_p f (x), \quad \text{for } \  l \geq 0 \ \text{ and } \ x \in [0,3^l-1].
\end{equation*}
\text{ } \\
Note that for $l=0$, the probabilistic graph Laplacian $\Delta^{(0)}_{p}$ is independent of the parameter $p$, and therefore we omit it from the notation in this case
\begin{equation} 
\label{eq:probabilisticLaplacianLevel0}
\Delta^{(0)}: =\Delta^{(0)}_{p}=
\begin{pmatrix}
1 & -1\\
-1 & 1
\end{pmatrix}.
\end{equation}
\begin{figure}[htp]
    \centering
  \includegraphics[width=1.\textwidth]{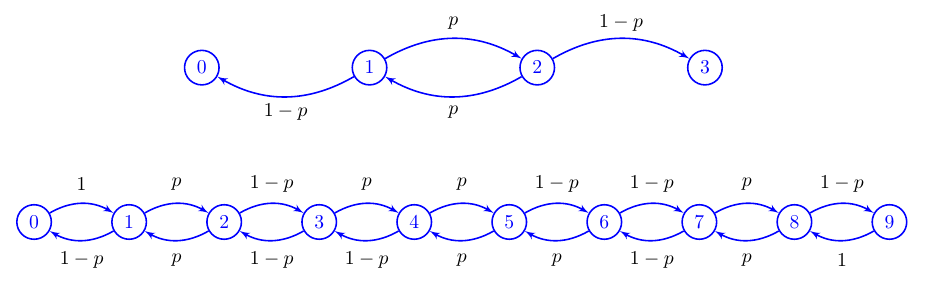}
\caption{ Visual illustration of the substitution rule. (Top) A copy of $G_1$. The deleted edges correspond to the edges that are replaced when applying the substitution rule. (Bottom) The graph $G_2$, which is constructed by inserting the three copies of $G_1$ in protograph shown in Figure \ref{fig:modelGraphNeumann}. While the vertices are labeled by the addresses, the labeling of the edges represents the transition probabilities (off-diagonal entries in the self-similar Laplacian).}
\label{fig:G2Neumann}
\end{figure}

\subsection{The self-similar almost Mathieu operators}\label{subsec2.2}

We introduce a \textit{self-similar version of almost Mathieu operators} defined with respect to the self-similar Laplacian $\Delta_p$ introduced in the last section. Let $f \in \ell(\mathbb{Z}_+)$, $\alpha \in \rr$, $\theta \in [0,2 \pi)$ and $\beta \in \rr$. We define
\begin{align}
	\label{eq:AMOversionPQ}
	(\hamilton_{p,\beta, \alpha, \theta } f)(x) = \left\{\begin{array}{ll} \beta \cos{( \theta)} f(0)-f(1), & \text{if}~x=0. \\
		\beta \cos{(2 \pi \alpha x + \theta)} f(x) & p(x,x-1)=1-p, \ p(x,x+1)=p, \\
		\quad \quad -p(x,x-1)f(x-1)-p(x,x+1)f(x+1), & \text{if }~3^{-m(x)}x \equiv 1\pmod 3.
		\\
		\beta \cos{(2 \pi \alpha x + \theta)} f(x) &  p(x,x-1)=p, \ p(x,x+1))=1-p, \\
		\quad \quad -p(x,x-1)f(x-1)-p(x,x+1)f(x+1),&\text{if}~3^{-m(x)}x \equiv 2\pmod 3.
	\end{array}\right..
\end{align}
Setting $p=\frac{1}{2}$ recovers up to a multiplicative constant the common form of the one-dimensional almost Mathieu operators, i.e. for $x \in \mathbb{Z}_+ \backslash \{0\}$, 
\begin{align}
\label{eq:UsualAMO}
(\hamilton_{\frac{1}{2},\beta, \alpha, \theta } f)(x) = -\frac{1}{2} \Big( f(x+1) + f(x-1) -2 \beta \cos{(2 \pi \alpha x + \theta)} f(x) \Big).
\end{align}
By Lemma \ref{lem:SelfAdjointnessOfLap}, $\hamilton_{p,\beta, \alpha, \theta}$  is a bounded self-adjoint operator on $ \ell^2(\integers_{+}, d \pi)$. 

For the sequence of graphs $\{G_l\}_{l \in \nn}$ given in Definition~\ref{def:finiteGraphApproxNeumann}, we associate a truncation $\hamilton_{p,\beta, \alpha, \theta }^{(l)} := \hamilton_{p,\beta, \alpha, \theta }|_{V_l} $, of the almost Mathieu operators~\eqref{eq:AMOversionPQ}, where we recall, $V_{l} = \mathbb{Z}_+ \cap [0,3^l]$. In particular, $\hamilton_{p,\beta, \alpha, \theta }^{(l)} $ is given by 

\begin{align}
\label{eq:truncatedAMO}
	(\hamilton_{p,\beta, \alpha, \theta }^{(l)} f)(x) = \left\{\begin{array}{ll}  \beta \cos{( \theta)} f(0)-f(1), & \text{if}~x=0. \\
		 \beta \cos{(2 \pi \alpha 3^l + \theta)} f(3^l)-f(3^l-1), & \text{if}~x=3^l. \\
		\beta \cos{(2 \pi \alpha x + \theta )} f(x) & p(x,x-1)=1-p, \ p(x,x+1)=p, \\
		\quad \quad -p(x,x-1)f(x-1)-p(x,x+1)f(x+1), & \text{if }~3^{-m(x)}x \equiv 1\pmod 3.
		\\
		\beta  \cos{(2 \pi \alpha x + \theta)} f(x) &  p(x,x-1)=p, \ p(x,x+1))=1-p, \\
		\quad \quad -p(x,x-1)f(x-1)-p(x,x+1)f(x+1),&\text{if}~3^{-m(x)}x \equiv 2\pmod 3.
	\end{array}\right..
\end{align}
\text{ } \\
Note that, similarly to the construction of the $\{G_{l}\}_{l \geq 0}$, we impose a reflecting boundary condition on the vertex $3^l$. The restriction of $\hamilton_{p,\beta, \alpha, \theta }^{(l)}$ to the interior vertices of $G_l$ is denoted by $\hamilton_{p,\beta, \alpha, \theta }^{(l), D}$, i.e.

\begin{align}
\label{eq:DiriAMO}
\hamilton_{p,\beta, \alpha, \theta }^{(l)} = 
\left( 
\begin{array}{c | c | c} 
\beta \cos{( \theta)} & \begin{array}{c c c} 
     -1 & 0 & \dots 
  \end{array} & 0 \\
 \hline 
\vdots &  \begin{array}{c c c} 
      &  &    \\ 
      & \DiriAMO & \\ 
      &  &  
  \end{array} & \vdots \\ 
  \hline 
0 &  \begin{array}{c c c} 
     \dots & 0 & -1
  \end{array} &  \beta \cos{(2 \pi \alpha 3^l + \theta)}
 \end{array} 
\right).
 \end{align}
 \text{ } \\
We identify $\hamilton_{p,\beta, \alpha, \theta }^{(l),D}$ with $\hamilton_{p,\beta, \alpha, \theta }^{(l)}$ when defined on
 the domain $\{ f:V_{l} \to \complex \ | \  f(0)=f(3^l)=0 \ \}$. We refer to $\hamilton_{p,\beta, \alpha, \theta }^{(l),D}$ as the \textit{Dirichlet almost Mathieu operator of level $l$}.
 In the following, we regard the matrix $\hamilton_{p,\beta, \alpha, \theta }^{(l)}$  as extended by zeros to a semi-infinite matrix. 
\begin{proposition}
Let $f \in \ell^2(\integers_{+}, d \pi) $. Then 
$$\begin{cases} \ \lim_{l \to \infty} ||   \hamilton_{p,\beta, \alpha, \theta } f -  \hamilton_{p,\beta, \alpha, \theta }^{(l)} f || = 0\\
 \ \lim_{l \to \infty} ||      \big( z-  \hamilton_{p,\beta, \alpha, \theta }  \big)^{-1} f -  \big( z-  \hamilton_{p,\beta, \alpha, \theta }^{(l)}  \big)^{-1} f || = 0.\end{cases}$$
\end{proposition}
The strong convergence is evident and strong resolvent convergence follows by \cite{Weidmann1997}. The reader is also referred to \cite{reed1981functional}. We note that the statement holds as well for $\Delta^{(l)}_{p}$ and $\Delta_{p}$.


\section{Spectral analysis of the self-similar almost Mathieu operators}
\label{sec-spectralDeci}
In this section we prove  our two main results. First, we consider the truncated self-similar AMO $\hamilton_{p,\beta, \alpha, \theta }^{(l)}$, and prove that their spectra can be determined using the spectral decimation method when the parameter $\alpha$ is restricted to the set $\{\tfrac{k}{3^n},\, k=0,1, \hdots, 3^n -1\}_{n=1}^l$ where $l\geq 1$ is the truncation level. In particular, this finite graph case is given in Theorem~\ref{thm:firstResultFiniteGraphs}. Subsequently, we state Theorem~\ref{thm:SpectralSimiAMOandpqModelnew} under the same restriction on the parameter $\alpha$.

\subsection{Finite graphs case}\label{sec-finitegraphs}
This section will briefly review a now standard technique  used in Analysis on Fractals and called \textit{Spectral Decimation}. We prove that it can be applied to  the sequence of almost Mathieu operators $\hamilton_{p,\beta, \alpha, \theta }^{(l)}$ for $\alpha = \frac{k}{3^n}$, $k \in \integers$, $1 \leq n \leq l$ and $\theta =0$. The method was intensively applied in the context of  Laplacians on fractals and self-similar graphs. Its central idea is that the spectrum of such Laplacian can be completely described in terms of iterations of a rational function,  called the \textit{spectral decimation function}. Below, we extend this method to the self-similar almost Mathieu operators when the frequency $\alpha$ is appropriately calibrated with the hierarchical structure of the self-similar Laplacian. In this case, we provide a complete description of the spectrum of $l$th-level almost Mathieu operators $\hamilton_{p,\beta, \frac{k}{3^n}, 0 }^{(l)}$ by relating it to the spectrum of $(l-n)$th-level Laplacian, i.e. $\sigma(\Delta^{(l-n)}_{p})$. The following theorem is the main result of this section.

 \begin{theorem}
\label{thm:firstResultFiniteGraphs}
Let $p\in (0,1)$, $\beta \in \rr$, $l \geq 1$, and $1 \leq n \leq l$ be fixed. Let $\theta =0$, and for $k \in \{0,\dots, 3^n-1\}$ set $\alpha = \frac{k}{3^n}$. There exists a polynomial $R_{p,\beta, \frac{k}{3^n},0 }$ of order $3^n$ such that,
\begin{equation}
\sigma
\Big( \hamilton_{p,\beta, \frac{k}{3^n}, 0 }^{(l)}  \Big)   = R^{-1}_{p,\beta, \frac{k}{3^n},0 } \Big( \sigma(\Delta^{(l-n)}_{p}) \backslash \sigma(\Delta^{(0)})  \Big) \bigcup  \sigma \Big( \hamilton_{p,\beta, \frac{k}{3^n}, 0 }^{(n)}  \Big).
\end{equation}
\text{ } \\
Furthermore, for $n=1$ and $k \in \{1,2\}$, the polynomial is given by

\begin{equation*}
R_{p,\beta, \frac{k}{3}, 0 }(z) =  \frac{\left(  - \beta + 2 p - 2 z  \right) \left(\beta^{2} + 2 \beta p + \beta z - 2 p z - 2 p - 2 z^{2} + 2\right)}{4 p \left(1-p\right)}.
\end{equation*}
 \end{theorem}
 
 Before giving the proof of this result, we recall some facts that can be found in~\cite{malozemovteplyaev2003}. Let $\mathcal{H}$ and $\mathcal{H}_0$ be Hilbert spaces, and $U:\mathcal{H}_0 \to \mathcal{H}$ be an isometry. Suppose  $H$ and $H_0$ are bounded linear operators on  $\mathcal{H}$ and $\mathcal{H}_0$, respectively, and that $\phi,\psi$ are complex-valued functions. Following \cite[Definition 2.1]{malozemovteplyaev2003}, we say that the operator $H$ \textit{spectrally similar} to the operator $H_0$ with functions $\phi$ and $\psi$ if 
 \begin{equation}
\label{eq:OriginalSpectSimi}
U^{\ast}(H-z)^{-1}U=(\phi(z)H_0 - \psi(z))^{-1},
\end{equation}
for all $z \in \complex$ such that the two sides of~\eqref{eq:OriginalSpectSimi} are well defined. In particular,  for $z$ in the domain of both $\phi$ and $\psi$ such that $\phi(z)\neq0$, we have $z\in\rho(H)$ (the resolvent of $H$) if and only if $R(z)=\frac{\psi(z)}{\phi(z)}\in\rho(H_0)$ (the resolvent of $H_0$).  We call $R(z)$ the {\em spectral decimation function}. In general, the functions $\phi(z)$ and $\psi(z)$ are usually difficult to express, but they can be computed effectively using the notion of  Schur complement. We refer to~\cite{malozemovteplyaev2003,BajorinSteinhurst2008,Bajorin_2007} for some examples. Identifying $\mathcal{H}_0$  with a closed subspace of $\mathcal{H}$ via $U$, let $\mathcal{H}_1$ be the orthogonal complement and decompose $H$ on   $\mathcal{H}=\mathcal{H}_0 \oplus \mathcal{H}_1$ in the block form
\begin{equation}
\label{eq:lDecompo}
H=\begin{pmatrix}
T & J^T\\
J & X
\end{pmatrix}.
\end{equation}
\begin{lemma}[\cite{malozemovteplyaev2003}, Lemma 3.3]
\label{lem:Lemma33}
For $z \in \rho(H)\cap \rho(X)$ the operators $H$ and $H_0$ are spectrally similar if and only if
the Schur complement of $H-zI$, given by $S_{H}(z)=T-z -  J^T(X-z)^{-1}J$, satisfies
\begin{equation}
\label{eq:schurComp1}
    S_{H}(z)=  \phi(z) H_0 - \psi(z)I.
\end{equation}
\end{lemma}

 $\mathscr{E}_{H}:=\{z \in \complex \ | \ z \in \sigma(X) \text{ or } \phi(z)=0 \}$ is called the \textit{exceptional set} of $H$ and  plays a crucial role in the spectral decimation method.  The spectral decimation has been already implemented for $\{ \Delta^{(n)}_{p}\}_{n \geq 0}$. For the sake of completeness we state this result and refer to \cite[Lemma 5.8]{TeplyaevSpectralZeta2007} for more details.

 \begin{proposition}\cite[Lemma 5.8]{TeplyaevSpectralZeta2007}
 \label{prop: SpectralDecimpLap}
 Let $n \geq 1$, then $\Delta^{(n)}_{p}$ is spectrally similar to $\Delta^{(n-1)}_{p}$ (with respect to functions given in \cite{TeplyaevSpectralZeta2007}). The spectral decimation function $R_{\Delta_{p}}$ and the exceptional set is $$\mathscr{E}_{\Delta_{p}}=\{1+p, 1-p \}$$ and 
\begin{equation}
\label{eq:specDeciFuncLap}
R_{\Delta_{p}}(z) =  \frac{z^3-3z^2+(2+p(1-p))z}{p(1-p)} =  \frac{z(z+p-2)(z-p-1)}{p(1-p)}
\end{equation}
Moreover, $\sigma(\Delta^{(0)})=\{0,2\}$ and  $\sigma(\Delta^{(n)}_{p})=\sigma(\Delta^{(0)}) \cup \bigcup_{i=0} ^{n-1} R_{\Delta_{p}}^{-i} (\{ p,2-p\}) $ for $n \geq 1$.
 \end{proposition}
 
 For the rest of this section, we fix  $p\in (0,1)$, $\beta \in \rr$, $l \geq 1$, and $1 \leq n \leq l$. We set $\theta =0$, $k \in \{0,\dots, 3^n-1\}$ and $\alpha = \frac{k}{3^n}$. We apply Lemma \ref{lem:Lemma33} on the level $l$ almost Mathieu operator $\hamilton_{p,\beta, \frac{k}{3^n}, 0 }^{(l)}$. We obtain the block form (\ref{eq:lDecompo}) by decomposing $\hamilton_{p,\beta, \frac{k}{3^n}, 0 }^{(l)}$ with respect to
\begin{align}
\mathcal{H}_0:=\text{span}\{\delta_v \ | \ v \text{ mod } 3^n \equiv 0 \} ,\quad \mathcal{H}_1:=\text{span}\{\delta_v \ | \ v \text{ mod } 3^n \not\equiv 0 \}.
\end{align}
where $\{\delta_x\}_{x \in V_{l}}$ is the canonical basis defined in (\ref{eq:canonicalBasis}) and $V_{l} = \mathbb{Z}_+ \cap [0,3^l]$. 
In practical terms:
\begin{enumerate}
	\item We rearrange the vertices in such a way that all vertices $v \in V_l$ with $v \text{ mod } 3^n \equiv 0$ appear before all vertices with $v \text{ mod } 3^n \not\equiv 0$, i.e. $V_l=\{0,3^n,\dots, 3^l, 1, 2, .. , 3^l-1\}$.
	\item We represent the matrix $\hamilton_{p,\beta, \frac{k}{3^n}, 0 }^{(l)}$ with respect to the canonical basis so that the order of the basis vectors follows the order of the vertices in step one.
	\item The matrix $\hamilton_{p,\beta, \frac{k}{3^n}, 0 }^{(l)}$ is then decomposed into the following block form
\begin{equation}
\label{eq:lDecompoLevelL}
\hamilton_{p,\beta, \frac{k}{3^n}, 0 }^{(l)}
=
\begin{pmatrix}
T_l & J_l^T\\
J_l & X_l
\end{pmatrix},
\end{equation}
where $T_l$ and $X_l$, correspond to the basis vectors $\{\delta_v \ | \ v \text{ mod } 3^n \equiv 0 \}$ and $\{\delta_v \ | \ v \text{ mod } 3^n \not\equiv 0 \}$, respectively.
\end{enumerate}
We observe that $T_l$ is a multiple of the identity matrix and that $X_l$ is a block diagonal matrix in which the diagonal blocks are the $n$th level Dirichlet almost Mathieu Operator $\hamilton_{p,\beta, \frac{k}{3^n}, 0 }^{(n),D}$, i.e.
 \begin{equation}
 \label{eq: decompoMatrices}
T = \beta 
\left(\begin{matrix}
1 & 0 &  \dots & 0 \\
0 & 1 &  \dots & 0 \\
\vdots & \vdots & \ddots & \vdots  \\
0 & 0 &  \dots  & 1
\end{matrix}\right), \quad \quad 
X_l =
  \renewcommand*{\arraystretch}{1.4}
 \left(
\begin{array}{@{}cc@{}c@{}c}
 \ \DiriAMOn
  & & & \\
   & \DiriAMOn & & \\
 &     & \  \bigdots \ \text{ } & \\
&&& \DiriAMOn
\end{array}
 \right).
 \end{equation}
In particular, we imply $\sigma(X_l)= \sigma(\hamilton_{p,\beta, \frac{k}{3^n}, 0 }^{(n),D})$.

\begin{lemma}
\label{lem:existOfphiandVerphi}
Let $p\in (0,1)$, $\beta \in \rr$, $l \geq 1$, and $1 \leq n \leq l$ be fixed. Moreover, we set $\theta =0$,  $\alpha = \frac{k}{3^n}$, for $k \in \{0,\dots, 3^n-1\}$. There exist functions  $ \phi_{p,\beta, \frac{k}{3^n},0}$ and $\psi_{p,\beta, \frac{k}{3^n},0 }$, such that $\hamilton_{p,\beta, \frac{k}{3^n}, 0 }^{(l)}$ is spectrally similar to $\Delta^{(l-n)}_{p}$ with respect to $ \phi_{p,\beta, \frac{k}{3^n},0}$  and $\psi_{p,\beta, \frac{k}{3^n},0}$.
\end{lemma}

\begin{proof}
Due to \cite[Lemma 3.10]{malozemovteplyaev2003}, it is sufficient to prove the existence of such functions  $ \phi_{p,\beta, \frac{k}{3^n},0} $ and $\psi_{p,\beta, \frac{k}{3^n},0 }$, so that the $n$th level $\hamilton_{p,\beta, \frac{k}{3^n} , 0 }^{(n)}$ is spectrally similar  to
\begin{equation} 
\label{eq:probabilisticLaplacianLevel0new}
\Delta^{(0)} =
\begin{pmatrix}
1 & -1\\
-1 & 1
\end{pmatrix}
\end{equation}
with the same functions $ \phi_{p,\beta, \frac{k}{3^n},0} $ and $\psi_{p,\beta, \frac{k}{3^n},0 }$.
The assumption  $\alpha = \frac{k}{3^n}$ guarantees that the matrix $\hamilton_{p,\beta, \frac{k}{3^n} , 0 }^{(n)}$ is symmetric with respect to its boundary vertices in the sense of  \cite[Definition 4.1]{malozemovteplyaev2003}. The spectral similarity of $\hamilton_{p,\beta, \frac{k}{3^n} , 0 }^{(n)}$ and $\Delta^{(0)}$ follows then by  \cite[Lemma 4.2]{malozemovteplyaev2003}.
\end{proof}

\begin{remark}
As a domain of   $ \phi_{p,\beta, \frac{k}{3^n},0} $ and $\psi_{p,\beta, \frac{k}{3^n},0 }$  we use the resolvent $\rho(X_l)$
of $X_l$, where $X_l$ is the block diagonal matrix in~\eqref{eq: decompoMatrices}. 
For more details about this facts we refer to~\cite[Corollary 3.4]{malozemovteplyaev2003}.
\end{remark}

\begin{proposition}
\label{prop:PropertiesOfR} Let $p\in (0,1)$, $\beta \in \rr$, $l \geq 1$, and $1 \leq n \leq l$ be fixed, and  set $\theta =0$,  $\alpha = \frac{k}{3^n}$, for $k \in \{0,\dots, 3^n-1\}$. The following statements hold:
\begin{enumerate}
	\item $\phi_{p,\beta, \frac{k}{3^n},0}(z) \neq 0 $ for all $z \in \rho(X_l)$.
	\item   The exceptional set of $\hamilton_{p,\beta, \frac{k}{3^n}, 0 }^{(l)}$ is given by $\mathscr{E}_{p,\beta, \frac{k}{3^n},0}=\sigma(\hamilton_{p,\beta, \frac{k}{3^n}, 0 }^{(n),D})$.
	\item The spectral decimation function $R_{p,\beta, \frac{k}{3^n},0}(z):=\frac{\psi_{p,\beta, \frac{k}{3^n},0 }(z)}{\phi_{p,\beta, \frac{k}{3^n},0}(z) }$ is a polynomial of order $3^n$.
	\item $z \in \sigma \big(\hamilton_{p,\beta, \frac{k}{3^n}, 0 }^{(n)}\big) \bigcup \sigma \big(\hamilton_{p,\beta, \frac{k}{3^n}, 0 }^{(n),D} \big)$ if and only if $R_{p,\beta, \frac{k}{3^n},0}(z) \in  \sigma(\Delta^{(0)})$.
\end{enumerate}
\end{proposition}

\begin{proof}
 We prove this result in a more general setting of mirror-symmetric Jacobi matrices in a companion paper \cite{BaluMograbyOkoudjouTeplyaevJacobi2021}.
\end{proof}

The following could be derived immediately from Lemma~\ref{lem:existOfphiandVerphi}, but for the sake of completeness and clarity we give the details leading to  explicit formulas for $\phi_{p,\beta,\frac{k}{3},0}$, $\psi_{p,\beta, \frac{k}{3} ,0}$ and $ R_{p,\beta,\frac{k}{3},0 }$.
\begin{lemma}
\label{prop:SpectralSimilarityToLap}
Let $n=1$ and $k \in \{1,2\}$. Then $\hamilton_{p,\beta, \frac{k}{3}, 0 }^{(l)}$ is spectrally similar to  $\Delta^{(l-1)}_{p}$ with the functions

\begin{equation}
\label{eq:phiAndpsi}
 \phi_{p,\beta,\frac{k}{3},0}(z)=\frac{4 p \left(p - 1\right)}{4 p^{2} - \left( \beta + 2 z\right)^{2}} , 
  \quad \quad \psi_{p,\beta, \frac{k}{3},0 }(z) =  - \frac{\beta^{2} + 2 \beta p + \beta z - 2 p z - 2 p - 2 z^{2} + 2}{\beta + 2 p + 2 z}.
\end{equation}
\text{ } \\
The spectral decimation function $R_{p,\beta,\frac{k}{3},0 }$ and the exceptional set $\mathscr{E}_{p,\beta, \frac{k}{3},0}$ are given by
\begin{equation}
\label{eq:SpecDeciandExcepSet}
 R_{p,\beta,\frac{k}{3},0 }(z) =  \frac{\left(  - \beta + 2 p - 2 z  \right) \left(\beta^{2} + 2 \beta p + \beta z - 2 p z - 2 p - 2 z^{2} + 2\right)}{4 p \left(1-p\right)} , 
  \quad \quad \mathscr{E}_{p,\beta, \frac{k}{3},0} = \left\{ - \frac{\beta}{2} - p, \  - \frac{\beta}{2} + p \right\}.
\end{equation}
 \end{lemma}

 \begin{proof}
With the same argument as in the proof of Lemma \ref{lem:existOfphiandVerphi}, it is sufficient to consider the spectral similarity between $\hamilton_{p,\beta, \frac{k}{3} , 0 }^{(1)}$ and $\Delta^{(0)}$. Applying the above three steps on the level-one almost Mathieu operator gives
 \begin{equation}
\label{eq:natricesLevel1and2}
\hamilton_{p,\beta, \frac{1}{3}, 0 }^{(1)} =
\renewcommand*{\arraystretch}{1.3}
\left(
\begin{array}{cc|cc}
\beta & 0 & -1 & 0 \\
0 & \beta & 0 & -1 \\
\hline
p - 1 & 0 & - \frac{ \beta}{2} & - p\\
0 & p - 1 & - p & - \frac{ \beta}{2}
\end{array}
\right)
,\quad
X_1=
    \left(\begin{array}{c c}
   - \frac{ \beta}{2} & - p \\
  - p & - \frac{ \beta}{2}
\end{array}\right).
\end{equation}
We compute the Schur complement and express it as a linear combination $   \phi_{p,\beta, \frac{k}{3},0}(z) \Delta^{(0)}- \phi_{p,\beta, \frac{k}{3},0}(z)I$, 

\begin{equation}
\label{eq:SchurMatrix}
\renewcommand*{\arraystretch}{1.3}
 \left(
\begin{array}{cc}
\beta - z + \frac{\left(\frac{ \beta}{2} + z\right) \left(p - 1\right)}{p^{2} - \left(\frac{ \beta}{2} + z\right)^{2}} & -\frac{4 p \left(p - 1\right)}{4 p^{2} - \left( \beta + 2 z\right)^{2}}\\    -\frac{4 p \left(p - 1\right)}{4 p^{2} - \left( \beta + 2 z\right)^{2}} & \beta - z + \frac{\left(\frac{ \beta}{2} + z\right) \left(p - 1\right)}{p^{2} - \left(\frac{ \beta}{2} + z\right)^{2}}
\end{array}
\right)=  
 \phi_{p,\beta, \frac{k}{3},0 }(z) \left(
\begin{array}{cc}
 1 & -1 \\
-1 & 1
\end{array}
\right)
-
\varphi_{p,\beta, \frac{k}{3},0 }(z) \left(
\begin{array}{cc}
1 & 0 \\
 0& 1
\end{array}
\right).
\end{equation}
\text{} \\
The formulas (\ref{eq:phiAndpsi}) and (\ref{eq:SpecDeciandExcepSet}) can be verified by comparing both sides of the equation (\ref{eq:SchurMatrix}).
\end{proof}

\begin{proof}[Proof of Theorem \ref{thm:firstResultFiniteGraphs}] We note that the spectra of $\{\Delta^{(n)}_{p}\}^{\infty}_{n=0}$ are nested, i.e. $\{0,2\} = \sigma(\Delta^{(0)}) \subset \sigma(\Delta^{(1)}_{p})  \subset \dots \subset [0,2]$.  We split the preimages set into two subsets:
\begin{enumerate}
	\item  $R^{-1}_{p,\beta, \frac{k}{3^n},0 } \Big( \sigma(\Delta^{(l-n)}_{p}) \backslash \sigma(\Delta^{(0)}) \Big)$: There are $3^{(l-n)}+1$ distinct eigenvalues in $\sigma(\Delta^{(l-n)}_{p})$. In particular, $\big| \sigma(\Delta^{(l-n)}_{p}) \backslash \sigma(\Delta^{(0)})\big|=3^{(l-n)}-1$ and 
	\begin{align*}
	\Big|R^{-1}_{p,\beta, \frac{k}{3^n},0 } \Big( \sigma(\Delta^{(l-n)}_{p}) \backslash \sigma(\Delta^{(0)}) \Big) \Big|=3^n(3^{(l-n)}-1) = 3^l-3^n.
	\end{align*}
	Note that by Proposition \ref{prop:PropertiesOfR}(4), we conclude that all the $3^l-3^n$ preimages are not in the exceptional set and therefore eigenvalues of  $\hamilton_{p,\beta, \frac{k}{3^n}, 0 }^{(l)}$, see \cite[Theorem 3.6.(2)]{malozemovteplyaev2003}. Besides, this implies that all the $3^l-3^n$ preimages are distinct eigenvalues.
	\item  $R^{-1}_{p,\beta, \frac{k}{3^n} ,0} \Big( \sigma(\Delta^{(0)}) \Big)$: By Proposition \ref{prop:PropertiesOfR}(4), we have
	\begin{equation*}
	 R^{-1}_{p,\beta, \frac{k}{3^n},0 } \big( \sigma(\Delta^{(0)}) \big)=\sigma \big(\hamilton_{p,\beta, \frac{k}{3^n}, 0 }^{(n)}\big) \bigcup \sigma \big(\hamilton_{p,\beta, \frac{k}{3^n}, 0 }^{(n),D} \big),
	\end{equation*}
	By excluding the exceptional points, we see that $R^{-1}_{p,\beta, \frac{k}{3^n},0 } \big( \sigma(\Delta^{(0)}) \big)$ generates $3^n+1$ distinct eigenvalues of  $\hamilton_{p,\beta, \frac{k}{3^n}, 0 }^{(l)}$, namely the eigenvalues in $\sigma \big(\hamilton_{p,\beta, \frac{k}{3^n}, 0 }^{(n)}\big)$.
\end{enumerate}
We generated in part one and two $3^l-3^n + 3^n+1 = 3^l+1$ distinct eigenvalues, which shows with a dimension argument that we completely determined the spectrum $\sigma
\Big( \hamilton_{p,\beta, \frac{k}{3^n}, 0 }^{(l)}  \Big)$.
\end{proof}


\subsection{Infinite graphs case}\label{sec:infinitegraphs}

We extend the statement of Theorem \ref{thm:firstResultFiniteGraphs} to infinite graphs.  We provide a complete description of the spectrum of the almost Mathieu operators $\hamilton_{p,\beta, \frac{k}{3^n}, 0 }$ by relating it to the  self-similar Laplacian's spectrum $\sigma(\Delta_{p})$. The following theorem is the main result. 

\begin{theorem}\label{thm:SpectralSimiAMOandpqModelnew}
Let  $\hamilton_{p,\beta, \alpha, \theta}$ and $\Delta_p$ be given as in  (\ref{eq:AMOversionPQ}) and (\ref{pLaplace}). Let $p\in (0,1)$, $\beta \in \rr$ and $n \geq 1$ be fixed. We set $\theta =0$,  $\alpha = \frac{k}{3^n}$, for $k \in \{1,\dots, 3^n-1\}$. There exists a polynomial $R_{p,\beta, \frac{k}{3^n},0 }$ of order $3^n$ such that,

\begin{equation}
\sigma
\Big( \hamilton_{p,\beta, \frac{k}{3^n}, 0 }  \Big)   = R^{-1}_{p,\beta, \frac{k}{3^n},0 } \Big( \sigma(\Delta_p) \Big). 
\end{equation}
\text{ } \\
Moreover, $\hamilton_{p,\beta, \frac{k}{3^n}, 0 }$ has purely singularly continuous spectrum if $p\neq \frac{1}{2}$.
\end{theorem}
\begin{proof} 
The first part of 
theorem \ref{thm:SpectralSimiAMOandpqModelnew} is a   consequence of  \cite[Lemma 3.10]{malozemovteplyaev2003}. We proceed as in the previous section and apply the spectral decimation method. We set  $\mathcal{H} = \ell^2(\integers_{+}, d \pi) $ and $\mathcal{H}_0 = \ell^2(3^n \integers_{+}, d \pi) $, $n \geq 1$. Strictly speaking, the  self-similar Laplacian $\Delta_p$ in Theorem \ref{thm:SpectralSimiAMOandpqModelnew} is defined on $ \ell^2(3^n\integers_{+}, d \pi)$. To understand this, we follow  \cite[page 125]{BellissardRenormalizationGroup1992} and introduce  a dilation operator 
\begin{equation}
D:   \ell^2(3^n\integers_{+}, d \pi)\to  \ell^2(\integers_{+}, d \pi), \quad \quad (Df)(x)=f(3^n x),
\end{equation}
and its co-isometric adjoint
\begin{equation}
D^{\ast}:  \ell^2(\integers_{+}, d \pi) \to \ell^2(3^n \integers_{+}, d \pi), \quad \quad  (D^{\ast}f)(3^n x)=f(x).
\end{equation}
Next, we define the operator $\tilde{\Delta}_p$ on $ \ell^2(3^n \integers_{+}, d \pi) $ to be $\tilde{\Delta}_p = D^{\ast}\Delta_p D$. According to  \cite{ChenTeplyPQmodel2016}, $\tilde{\Delta}_p$ on $\ell^2(3^n \integers_{+}, d \pi)$ is isometrically equivalent to $\Delta_p$ on $ \ell^2(\integers_{+}, d \pi) $ and $\sigma(\tilde{\Delta}_p)=\sigma(\Delta_p)$. In the following, we will omit the tilde and refer to $\tilde{\Delta}_p$ by $\Delta_p$. We regard $\mathcal{H}_0 = \ell^2(3^n \integers_{+}, d \pi) $ as a subspace of  $\ell^2(\integers_{+}, d \pi)$ and introduce  
$\mathcal{H}_1$ as the orthogonal complement of $\mathcal{H}_0$ in $\mathcal{H}$. Then $\hamilton_{p,\beta, \frac{k}{3^n}, 0 }$ is decomposed with respect to $\mathcal{H}_0\oplus \mathcal{H}_1$ into the following block form
\begin{equation}
\label{eq:lDecompoInf}
\hamilton_{p, \beta, \alpha, \theta}=\begin{pmatrix}
T & J^T\\
J & X
\end{pmatrix}.
\end{equation}
\text{ } \\
We observe that $T$ is a multiple of the identity and that $X$ is a block diagonal semi-finite matrix in which the diagonal blocks are the $n$th level Dirichlet almost Mathieu Operator $\hamilton_{p,\beta, \frac{k}{3^n}, 0 }^{(n),D}$, i.e.

 \begin{equation}
 \label{eq: decompoMatricesInf}
T = \beta
\left(\begin{matrix}
1 & 0 & 0 & \dots \\
0 & 1 & 0 & \dots \\
0 & 0 & 1 & \ddots  \\
\vdots & \vdots & \vdots & \ddots  
\end{matrix}\right), \quad \quad 
X =
  \renewcommand*{\arraystretch}{1.4}
 \left(
\begin{array}{@{}cc@{}c@{}}
 \ \DiriAMOn
  & & \\
   & \DiriAMOn & \\
 &     & \  \bigdots \ \text{ }

\end{array}
 \right).
 \end{equation}

Similar to the proof of Lemma \ref{lem:existOfphiandVerphi}, the spectral similarity of $\hamilton_{p,\beta, \frac{k}{3^n} , 0 }^{(n)}$ and $\Delta^{(0)}$ implies the spectral similarity of $ \hamilton_{p,\beta, \frac{k}{3^n}, 0 }$ and $\Delta_p$ with the same $ \phi_{p,\beta, \frac{k}{3^n},0} $, $\psi_{p,\beta, \frac{k}{3^n},0 }$, $\mathscr{E}_{p,\beta, \frac{k}{3^n},0}$ and $R_{p,\beta, \frac{k}{3^n},0 }$. By \cite[Theorem 3.6]{malozemovteplyaev2003}, we see that  for $z \notin \mathscr{E}_{p,\beta, \frac{k}{3^n},0}$,
\begin{align*}
z \in \sigma\Big( \hamilton_{p,\beta, \frac{k}{3^n}, 0 }  \Big) \quad \Leftrightarrow \quad R_{p,\beta, \frac{k}{3^n},0 }(z) \in \sigma(\Delta_p) \quad \Leftrightarrow \quad  \quad z \in R^{-1}_{p,\beta, \frac{k}{3^n},0 } \big( \sigma(\Delta_p) \big).
\end{align*}
Next, we show $ \mathscr{E}_{p,\beta, \frac{k}{3^n},0} \subset \sigma\big( \hamilton_{p,\beta, \frac{k}{3^n}, 0 }  \big)$.  To this end we use Proposition \ref{prop:PropertiesOfR} (4), that is $ \mathscr{E}_{p,\beta, \frac{k}{3^n},0} \subset R^{-1}_{p,\beta, \frac{k}{3^n},0 } (0,2)$ and the fact that $0$ and $2$ are not isolated points in the spectrum $\sigma(\Delta_p)$. Let $z \in \mathscr{E}_{p,\beta, \frac{k}{3^n},0} \cap R^{-1}_{p,\beta, \frac{k}{3^n} ,0} (0)$. By a continuity argument, we can find a sequence $\{\lambda_m\}_{m \in \nn} \subset  \sigma(\Delta_p)$, $0 < \lambda_m <2 $, $\lambda_m \to 0$ and a partial inverse of
$R_{p,\beta, \frac{k}{3^n} ,0}$ (which we will denote by $R^{-1}_{p,\beta, \frac{k}{3^n} ,0}$ to avoid extra notation), such that
$R^{-1}_{p,\beta, \frac{k}{3^n},0 }(\lambda_m) \to z$.  Again with proposition \ref{prop:PropertiesOfR}(4) we have $R^{-1}_{p,\beta, \frac{k}{3^n},0 }(\lambda_m) \notin \mathscr{E}_{p,\beta, \frac{k}{3^n},0}$ for all $m \in \nn$ and imply by \cite[Theorem 3.6]{malozemovteplyaev2003} that
\begin{equation}
R^{-1}_{p,\beta, \frac{k}{3^n} ,0}(\lambda_m) \in   \sigma\Big( \hamilton_{p,\beta, \frac{k}{3^n}, 0 }  \Big)  \quad \forall \ m \in \nn.
\end{equation}
By closedness of the spectrum, we conclude that $z \in  \sigma\Big( \hamilton_{p,\beta, \frac{k}{3^n}, 0 }  \Big) $, see also Remark \ref{rem:continuityOfPreimages}. The same argument holds for $z \in \mathscr{E}_{p,\beta, \frac{k}{3^n},0} \cap R^{-1}_{p,\beta, \frac{k}{3^n},0 } (2)$. The second part of the statement followes by \cite[Theorem 1]{ChenTeplyPQmodel2016} combined with \cite[Theorem 3.6]{malozemovteplyaev2003}.
\end{proof}

\section{Integrated density of states}
\label{IDSsection}

\begin{figure}[!htb]
  \begin{minipage}[b]{0.5\linewidth}
    \centering
    \includegraphics[width=1.09\linewidth]{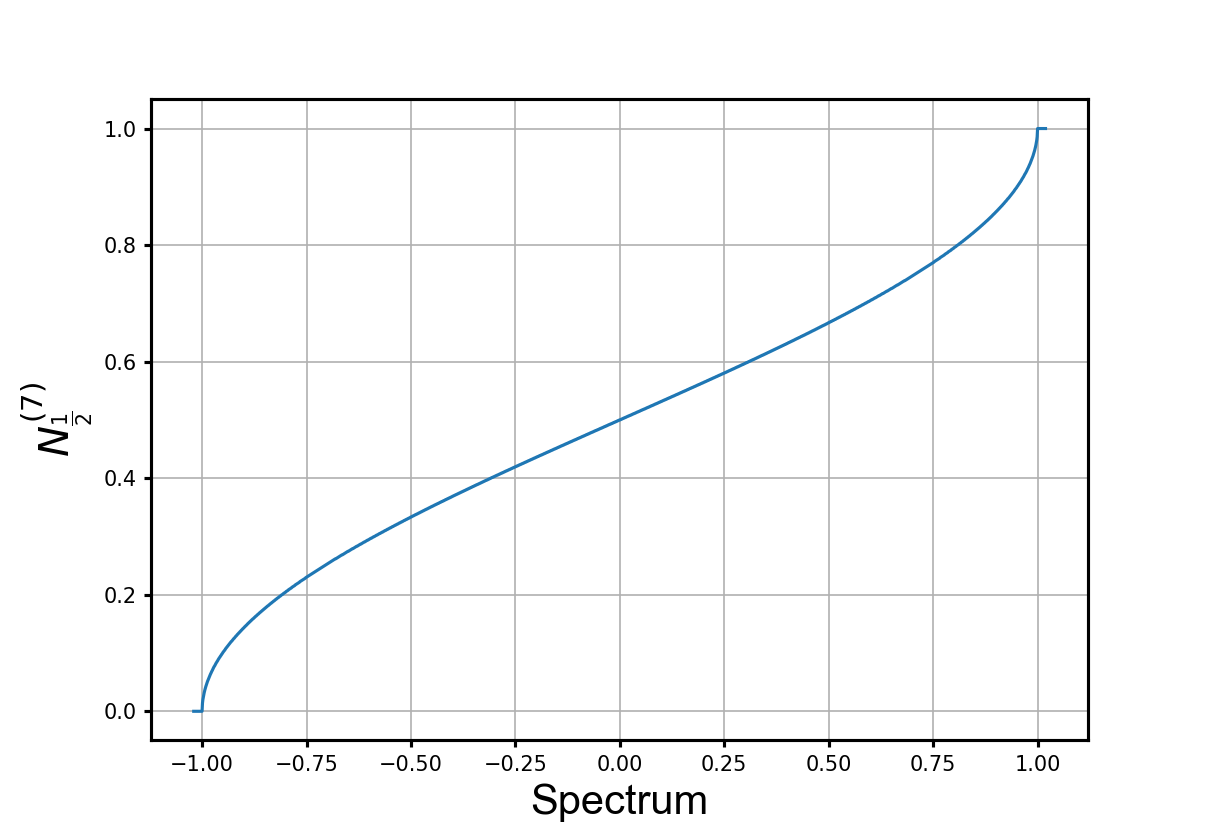} 
    \vspace{4ex}
  \end{minipage}
  \begin{minipage}[b]{0.5\linewidth}
    \centering
    \includegraphics[width=1.09\linewidth]{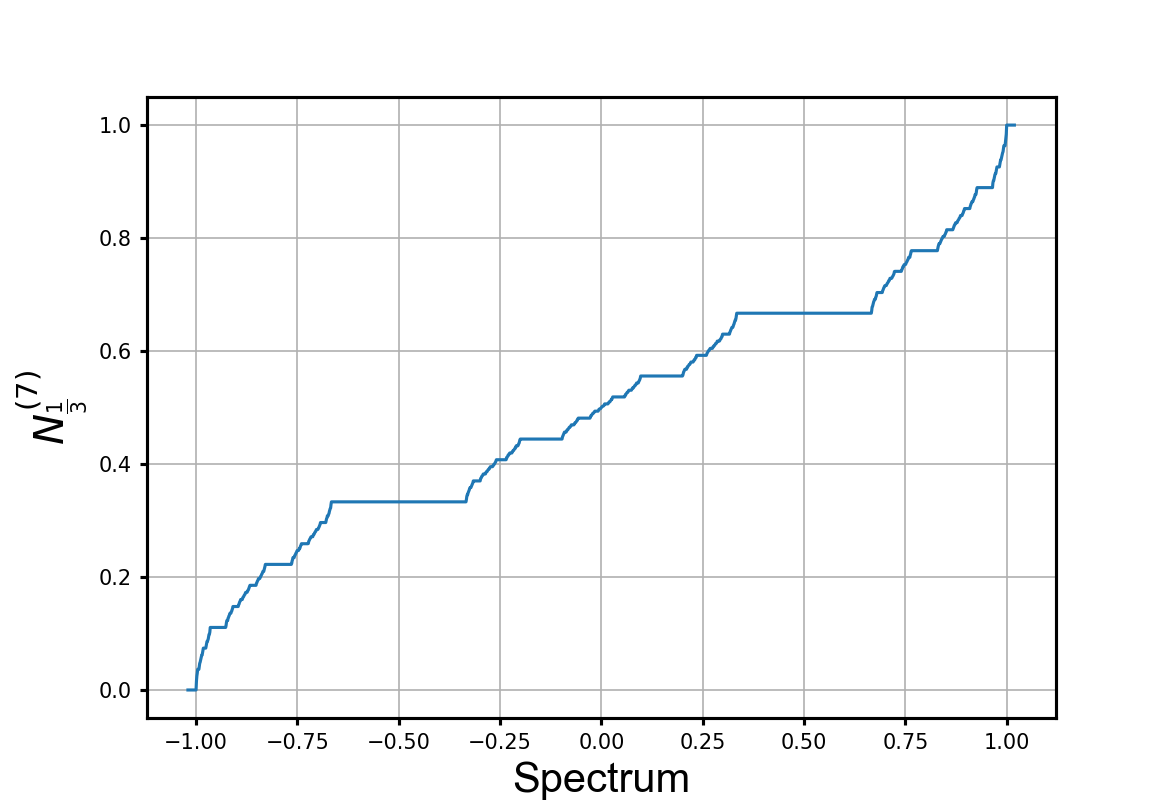} 
    \vspace{4ex}
  \end{minipage} 
  \caption{Numerical computation of the normalized eigenvalue counting function $N^{(l)}_{p}$. The computations are done for level $l=7$. (Left) $p=\frac{1}{2}$, i.e., standard probabilistic graph Laplacian with  $\sigma(\Delta_{\frac{1}{2}})=[0,2]$. (Right) $p=\frac{1}{3}$, i.e., a self-similar graph Laplacian where  $\sigma(\Delta_{\frac{1}{3}})$  is a cantor set.} 
  \label{fig:IDSforLaps}
\end{figure}

Throughout this section, we assume that $p\in (0,1)$ and $l \geq 1$  are fixed. We follow ideas presented in \cite[Section 5.4]{Kirsch2008Random} and define the \textit{density of states} of $\Delta_{p}$. We start by considering the spectrum of $\Delta^{(l)}_{p}$, which consists of finitely many simple eigenvalues. We refer to the normalized sum of Dirac measures concentrated on the eigenvalues
\begin{equation}
\nu_{l,p}(\{x\})
 =  \frac1{3^{l} +1} \sum_{\lambda\in  \sigma(\Delta^{(l)}_{p})} \delta_\lambda(x)
 \end{equation}
as the \textit{density of states} of $\Delta^{(l)}_{p}$. The \textit{normalized eigenvalue counting} function of  $\Delta^{(l)}_{p}$ is then given by $N^{(l)}_{p}(x):=\nu_{l,p}((-\infty,x])$.  We note that $\Delta^{(l)}_{p}$ is the restriction of $\Delta_{p}$ to the finite graph  $G_{l}=(V_{l}, E_{l})$ while imposing Neumann boundary conditions.  As the following results can be derived in the same way when Dirichlet boundary conditions are applied, we restrict our consideration to the former one.  Figure  \ref{fig:IDSforLaps} depicts the normalized eigenvalue counting function $N^{(l)}_{p}$ for different  parameters.
We recall some known facts about the spectrum of the self-similar Laplacian $\Delta_{p}$. Theorem 1 and Proposition 10  in \cite{ChenTeplyPQmodel2016} show that the spectrum $\sigma(\Delta_{p})$ is the Julia set $\mathcal{J}(R_{\Delta_{p}})$ of the polynomial $R_{\Delta_{p}}$ given in Proposition  \ref{prop: SpectralDecimpLap}, for more general settings see \cite{HareTep2011}. For $p=\frac{1}{2}$, we have $\mathcal{J}(R_{\Delta_{\frac{1}{2}}})=[0,2]$ and the spectrum is absolutely continuous. For $p \neq \frac{1}{2}$,  the Julia set $\mathcal{J}(R_{\Delta_{p}})$ is a Cantor set of Lebesgue measure zero and the spectrum is purely singularly continuous. 
Brolin in \cite{Brolin1965} proved the existence of a natural measure on polynomial Julia sets, namely the so-called balanced invariant measure. Moreover, he showed that the balanced invariant measure coincides with the potential theory's equilibrium (harmonic) measure. In higher generality the uniqueness of the balanced invariant measure was established later  in \cite{Freire1983,Ma1983},  the reader is referred to \cite{SmirnovThesis} for an overview. Denoting the balanced invariant measure of the Julia set $\mathcal{J}(R_{\Delta_{p}})$ by $\nu_p$ and  and using  ideas similar to Brolin's lead to  the following result.
\begin{proposition}
\label{prop:IDSofpLap}
The sequence of density of states $\{\nu_{l,p}\}_{l \in \nn}$ converges weakly to the balanced invariant measure $\nu_p$ of the Julia set $\mathcal{J}(R_{\Delta_{p}})$.
\end{proposition}
Let  $\theta =0$, $\beta \in \rr$,  $1 \leq n \leq l$ and $k \in \{0,\dots, 3^n-1\}$ be fixed. In the same way as above, we define the density of states and the normalized eigenvalue counting function of $\hamilton_{p,\beta, \frac{k}{3^n}, 0 }^{(l)}$ and refer to them by $\nu^{(l)}_{p,\beta, \frac{k}{3^n}, 0 }$ and $N^{(l)}_{p,\beta, \frac{k}{3^n}, 0 }$, respectively. Theorem \ref{thm:firstResultFiniteGraphs} asserts the existence of a polynomial $R_{p,\beta, \frac{k}{3^n},0 }$ of order $3^n$. Let $S_1,S_2, \dots, S_{3^n}$ be the $3^n$ branches of the inverse  $R^{-1}_{p,\beta, \frac{k}{3^n},0 }$ and $E \subset \rr$ ($ \nu_{p}$-measurable). We define the following measure
\begin{align}
\label{eq:DefNuForAMO}
\nu_{p,\beta, \frac{k}{3^n}, 0 }(E) :=\frac{1 }{3^{n} } \sum_{i=1}^{3^n} \int_{\sigma(\Delta_{p})} \chi_E(S_i(x)) \nu_{p}(dx)
\end{align}
where $\chi_E$ is the characteristic function of the set $E$, i.e.
\begin{equation}
\chi_E(x)
 =
  \begin{cases}
  \  0     & \quad \text{if } x \notin E \\ 
   \ 1     & \quad \text{if } x \in E
  \end{cases}
\end{equation}
\begin{theorem}
\label{thm:IDSofAMO}
Let $\text{supp}(\nu_{p,\beta, \frac{k}{3^n}, 0 })$ denotes the support of $\nu_{p,\beta, \frac{k}{3^n}, 0 }$. Then $\text{supp}(\nu_{p,\beta, \frac{k}{3^n}, 0 })=\sigma\big( \hamilton_{p,\beta, \frac{k}{3^n}, 0 }  \big) $. 
The sequence of density of states $\big\{\nu^{(l)}_{p,\beta, \frac{k}{3^n}, 0 } \big\}_{l \in \nn}$ converges weakly to $\nu_{p,\beta, \frac{k}{3^n}, 0 }$. Moreover, the following identity holds
\begin{align}
\label{eq:almostSelfSimiId}
\int_{\sigma( \hamilton_{p,\beta, \frac{k}{3^n}, 0 }  )} f(x) \nu_{p,\beta, \frac{k}{3^n}, 0 }(dx)=\frac{1 }{3^{n} } \sum_{i=1}^{3^n} \int_{\sigma(\Delta_{p})} f(S_i(x)) \nu_{p}(dx)
\end{align}
for  $f \in C_b(\complex)$, i.e. $f$ is a  continuous bounded function on $\complex$.
\end{theorem}

\begin{proof}
Let $f \in C_b(\complex)$. Theorem \ref{thm:firstResultFiniteGraphs} implies
\begin{align}
\label{eq:discIntAMOmeasure}
  \sum_{x \in \sigma( \hamilton_{p,\beta, \frac{k}{3^n}, 0 }^{(l)}) }  f(x) \nu^{(l)}_{p,\beta, \frac{k}{3^n}, 0 }(\{x\}) =  \frac{1}{3^l +1} \sum_{x \in R^{-1}_{p,\beta, \frac{k}{3^n},0 } ( \sigma(\Delta^{(l-n)}_{p}) \backslash \sigma(\Delta^{(0)})  ) }  f(x) +  \frac{1}{3^l +1} \sum_{x \in \sigma( \hamilton_{p,\beta, \frac{k}{3^n}, 0 }^{(n)}) }  f(x).
\end{align}
We show that the first term on the right-hand side of equation (\ref{eq:discIntAMOmeasure}) converges to the term on the right-hand side of equation (\ref{eq:almostSelfSimiId}).
\begin{align*}
\frac{1}{3^l +1} \sum_{x \in R^{-1}_{p,\beta, \frac{k}{3^n},0 } ( \sigma(\Delta^{(l-n)}_{p}) \backslash \sigma(\Delta^{(0)})  ) }  f(x) &  = \frac{1}{3^l +1} \sum_{x \in   \sigma(\Delta^{(l-n)}_{p}) \backslash \sigma(\Delta^{(0)})   } \sum_{i=1}^{3^n} f(S_i(x)) \\
&  = \frac{3^{l-n} +1}{3^l +1} \sum_{x \in   \sigma(\Delta^{(l-n)}_{p}) \backslash \sigma(\Delta^{(0)})   } \sum_{i=1}^{3^n} f(S_i(x)) \nu_{l-n,p}(\{x\})
\\
&  \to \frac{1 }{3^{n} } \sum_{i=1}^{3^n} \int_{\sigma(\Delta_{p})} f(S_i(x)) \nu_{p}(dx), \ \text{ as } \ l \to \infty.
\end{align*}
\end{proof}
\begin{remark}
\label{rem:continuityOfPreimages}
The existence and continuity of the branches of the inverse  $R^{-1}_{p,\beta, \frac{k}{3^n},0 }$ on the interval $[0,2]$ is given in a forthcoming work \cite{BaluMograbyOkoudjouTeplyaevJacobi2021}, where we develop a general framework by extending the results obtained in this paper to a large class of Jacobi operators. 
\end{remark}
Theorem \ref{thm:IDSofAMO} and proposition \ref{prop:IDSofpLap} justify the following definitions.
\begin{definition}
We refer to $\nu_{p,\beta, \frac{k}{3^n}, 0 }$ and $\nu_p$ as the \textit{density of states} of $\hamilton_{p,\beta, \frac{k}{3^n}, 0 }$ and $\Delta_{p}$, respectively. The corresponding \text{integrated density of states} are given by $N_{p,\beta, \frac{k}{3^n}, 0}(x)=\nu_{p,\beta, \frac{k}{3^n}, 0}((-\infty,x])$ and $N_{p}(x)=\nu_{p}((-\infty,x])$.
\end{definition}

It now follows that 

\begin{corollary}
\label{cor:scalingCopies}
Let $E \subset \sigma(\Delta_{p})$, then $\nu_{p,\beta, \frac{k}{3^n}, 0 }(S_j(E)) = \frac{1}{3^n} \nu_{p}(E)$.
\end{corollary}
\begin{proof}
We compute
\begin{align}
\nu_{p,\beta, \frac{k}{3^n}, 0 }(S_j(E)) =\frac{1 }{3^{n} } \sum_{i=1}^{3^n} \int_{\sigma(\Delta_{p})} \chi_{S_j(E)}(S_i(x)) \nu_{p}(dx)  =\frac{1 }{3^{n} }  \int_{\sigma(\Delta_{p})} \chi_{E}(x) \nu_{p}(dx),
\end{align}
where in the second equality, we use that $\chi_{S_j(E)}(S_i(x)) = 0$ ($\nu_p$ almost surely), whenever $i \neq j$ and $\chi_{S_j(E)}(S_j(x)) = \chi_{E}(x)$.
\end{proof}
The intuitively,   Theorem \ref{thm:IDSofAMO} implies that the density of states $\nu_{p,\beta, \frac{k}{3^n}, 0 }$ equally distributes the original mass $\nu_{p}$ of the spectrum $\sigma(\Delta_{p})$ on the $3^n$ branches of the inverse spectral decimation function $R^{-1}_{p,\beta, \frac{k}{3^n},0 }$. In particular, this enables us to compute the spectral gap labels of $\hamilton_{p,\beta, \frac{k}{3^n}, 0 }$. Let $\rho(\Delta_{p})$ be the resolvent set of $\Delta_{p}$. We define the set of spectral gap labels of $\Delta_{p}$ by
 \begin{align}
\mathscr{G} \mathscr{L}(\Delta_{p}) = \{ N_p(x) \ | \ x \in  \rho(\Delta_{p}) \cap \rr  \}.
\end{align}
It is not difficult to see that $\mathscr{G} \mathscr{L}(\Delta_{\frac{1}{2}}) =\{0,1\}$. For $p \neq \frac{1}{2}$, we have
 \begin{align}
\mathscr{G} \mathscr{L}(\Delta_{p}) = \Big\{ \frac{j}{3^i} \ \Big| \ i \in  \nn, \ j \in \{0,1,\dots,3^i\}  \Big\}.
\end{align}
We define the set of spectral gap labels of $\hamilton_{p,\beta, \frac{k}{3^n}, 0 }$ similarly and denote it by $\mathscr{G} \mathscr{L}(\hamilton_{p,\beta, \frac{k}{3^n}, 0 })$.
\begin{corollary}[Gap labeling]
\label{coro:GapLabel}
The set of spectral gap labels of $\hamilton_{p,\beta, \frac{k}{3^n}, 0 }$ is given by
\begin{align}
\mathscr{G} \mathscr{L}(\hamilton_{p,\beta, \frac{k}{3^n}, 0 }) \subset \Big\{ \frac{j}{3^n}+ \frac{1}{3^n}\mathscr{G} \mathscr{L}(\Delta_{p}) \ \Big| \  \ j \in \{0,1,\dots,3^n-1\}  \Big\}.
\end{align}
\end{corollary}

\section{Examples and numerical results}\label{sec:examplesappl}
\subsection{Spectra of  $\hamilton_{\frac{1}{3},1, \frac{1}{3}, 0 }^{(1)}$ and $\hamilton_{\frac{1}{3},1, \frac{1}{3}, 0 }^{(2)}$}
\begin{figure}[!htb]
    \centering
  \includegraphics[width=1.\textwidth]{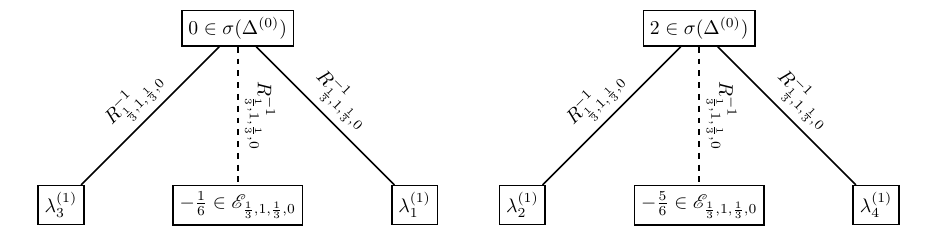}
    	\caption{The preimage sets $R^{-1}_{\frac{1}{3},1, \frac{1}{3}, 0 }(0) $ and $R^{-1}_{\frac{1}{3},1, \frac{1}{3}, 0 }(2) $. Note that $-\frac{1}{6}$ and $-\frac{5}{6}$ are elements of the exceptional set.   The numerical values are given in Table  \hyperref[tab:table1]{1}.}
	\label{fig:firstOffsprings}
\end{figure}

\begin{figure}[!htb]
    \centering
  \includegraphics[width=1.\textwidth]{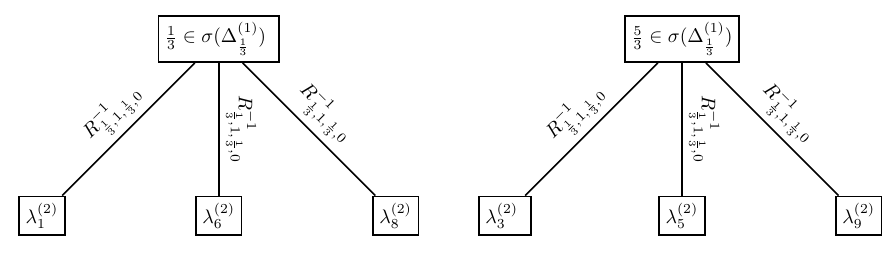}
    	\caption{The preimage sets $R^{-1}_{\frac{1}{3},1, \frac{1}{3}, 0 }(\frac{1}{3}) $ and $R^{-1}_{\frac{1}{3},1, \frac{1}{3}, 0 }(\frac{5}{3}) $.  The numerical values are given in Table \hyperref[tab: spectraForAlpha1over3first3levelsh2]{2}. }
	\label{fig:SecondOffsprings}
\end{figure}

We apply the above framework for finite graphs in the case $p=\frac{1}{3}$, $\beta=1$ and $\alpha = \frac{1}{3}$. Direct computations give $\sigma(\Delta^{(0)})=\{0,2\}$. With Proposition \ref{prop:SpectralSimilarityToLap} we compute the exceptional set and $R_{\frac{1}{3},1, \frac{1}{3}, 0} $,
\begin{equation}
\label{eq:spectralDecimationfunctionAndExceptset}
 \mathscr{E}_{\frac{1}{3},1, \frac{1}{3}, 0 } = \left\{ - \frac{1}{6} , \  - \frac{5}{6} \right\}, \quad \quad R_{\frac{1}{3},1, \frac{1}{3}, 0 }(z)= \frac{9 z^{3}}{2} - \frac{55 z}{8} - \frac{9}{8}.
\end{equation}
We give an illustration of Theorem~\ref{thm:firstResultFiniteGraphs}. Due to the spectral similarity between $\hamilton_{\frac{1}{3},1, \frac{1}{3}, 0 }^{(1)}$ and $\Delta^{(0)}$, we see that $z\in\sigma(\hamilton_{\frac{1}{3},1, \frac{1}{3}, 0 }^{(1)})\setminus\mathscr{E}_{\frac{1}{3},1, \frac{1}{3}, 0 }$ if and only if $R_{\frac{1}{3},1, \frac{1}{3}, 0 }(z)\in\sigma(\Delta^{(0)})$. We compute the preimage sets $R^{-1}_{\frac{1}{3},1, \frac{1}{3}, 0 }(0) $ and $R^{-1}_{\frac{1}{3},1, \frac{1}{3}, 0 }(2) $, see Figure  \ref{fig:firstOffsprings} and Table  \hyperref[tab:table1]{1}. We note that $R_{\frac{1}{3},1, \frac{1}{3}, 0 }$ is a polynomial of degree 3; therefore, each of the eigenvalues $0, 2 \in \sigma(\Delta^{(0)})$ generates three preimages. Excluding the exceptional points results in four distinct eigenvalues of  $\hamilton_{\frac{1}{3},1, \frac{1}{3}, 0 }^{(1)}$,  which on the other hand, determine the complete spectrum as $\hamilton_{\frac{1}{3},1, \frac{1}{3}, 0 }^{(1)}$ is a $4\times4$ matrix.
 \begin{table}[!htb]
\begin{tabular}{SSSSSS} \toprule
    {$\sigma(  \Delta^{(0)})$} & {$\lambda^{(0)}_1 = 0$} & {$\lambda^{(0)}_2 = 2$} &  &  & 
     \\ \midrule
    {$\sigma( \hamilton_{\frac{1}{3},1, \frac{1}{3}, 0 }^{(1)})$}  & {$\lambda^{(1)}_1 = \frac{1}{12} - \frac{\sqrt{217}}{12}$} & {$\lambda^{(1)}_2 =\frac{5}{12} - \frac{\sqrt{145}}{12}$} & {$\lambda^{(1)}_3 = \frac{1}{12} + \frac{\sqrt{217}}{12}$}  & {$\lambda^{(1)}_4 =    \frac{5}{12} + \frac{\sqrt{145}}{12}        $}   &   
     \\ \bottomrule \\
\end{tabular}
\label{tab:table1}
\caption{Numerical computation of the spectra $\sigma(  \Delta^{(0)})$ and $\sigma( \hamilton_{\frac{1}{3},1, \frac{1}{3}, 0 }^{(1)})$. The spectrum  $\sigma( \hamilton_{\frac{1}{3},1, \frac{1}{3}, 0 }^{(1)})$ is computed using Proposition \ref{prop:SpectralSimilarityToLap} and $\sigma(  \Delta^{(0)})$.}
\end{table} 
To compute $\sigma(\hamilton_{\frac{1}{3},1, \frac{1}{3}, 0 }^{(2)})$, we first use Proposition \ref{prop: SpectralDecimpLap}  and  the spectral decimation function $R_{\Delta_{p}}$ to calculate $\sigma(\Delta^{(1)}_{p})$. It can be easily checked that $\sigma(\Delta^{(1)}_{p})=\{0,\frac{1}{3},\frac{5}{3},2\}$. In particular, four out of the ten eigenvalues in $\sigma( \hamilton_{\frac{1}{3},1, \frac{1}{3}, 0 }^{(2)})$ are computed similarly to above, namely as the elements of preimage sets  $R^{-1}_{\frac{1}{3},1, \frac{1}{3}, 0 }(0) $ and $R^{-1}_{\frac{1}{3},1, \frac{1}{3}, 0 }(2) $ with excluding the points in the exceptional set. The preimage sets $R^{-1}_{\frac{1}{3},1, \frac{1}{3}, 0 }(1/3) $ and $R^{-1}_{\frac{1}{3},1, \frac{1}{3}, 0 }(5/3) $ are computed as shown in Figure  \ref{fig:SecondOffsprings} with the numerical values in Table \hyperref[tab: spectraForAlpha1over3first3levelsh2]{2}. These sets generate the remaining 6 eigenvalues. Note in level two, the graph $G_2$  consists of $10$ vertices.
\\ \\
\begin{table}[!htb]
\begin{tabular}{SSSSSS} \toprule
    {$\sigma(  \Delta^{(1)}_{\frac{1}{3}})$} & {$\lambda^{(1)}_1 = 0$} & {$\lambda^{(1)}_2 =\frac{1}{3}$} & {$\lambda^{(1)}_3 =\frac{5}{3}$} & {$\lambda^{(1)}_4 = 2$} &     \\ \midrule
    {$\sigma( \hamilton_{\frac{1}{3},1, \frac{1}{3}, 0 }^{(2)})$}  & {$\lambda^{(2)}_1 =-1.14424$} & {$\lambda^{(2)}_2 = -1.11189$} & {$\lambda^{(2)}_3 = -0.92631$}  & {$\lambda^{(2)}_4 = -0.58679$}   &    \\
    & {$\lambda^{(2)}_5 = -0.47717$}  & {$\lambda^{(2)}_6 =-0.21899$} & {$\lambda^{(2)}_7 = 1.31091$}  &  {$\lambda^{(2)}_8 = 1.33089$}  &
      \\
    & {$\lambda^{(2)}_9 = 1.40349$}  & {$\lambda^{(2)}_{10} =1.42013$} &   &   &
     \\ \bottomrule \\
\end{tabular}
\label{tab: spectraForAlpha1over3first3levelsh2}
\caption{ Numerical computation of the spectra $\sigma(  \Delta^{(1)}_{\frac{1}{3}})$ and $\sigma( \hamilton_{\frac{1}{3},1, \frac{1}{3}, 0 }^{(2)})$.}.
\end{table}

\subsection{Spectral gaps}\label{sec:spectralgaps}

The disconnectedness of the Julia set $\mathcal{J}(R_{\Delta_{p}})$, for $p \neq \frac{1}{2}$, implies that the self-similar Laplacian $\Delta_{p}$ has infinitely many spectral gaps. This fact combined with Theorem \ref{thm:SpectralSimiAMOandpqModelnew} lead us to the following two conclusions: \newline
\begin{enumerate}
\item For $p \neq \frac{1}{2}$, the spectrum $\sigma \Big( \hamilton_{p,\beta, \frac{k}{3^n}, 0 }  \Big)$ has infinitely many spectral gaps. 
\item We can generate the spectral gaps iteratively using the spectral decimation function $ R_{p,\beta, \frac{k}{3^n},0 }$. 
\end{enumerate}
We illustrate these ideas with the example $p=\frac{1}{3}$, $\beta=1$, $\alpha = \frac{1}{3}$, $\theta=0$ and generate the spectral gaps in $\sigma \big( \hamilton_{\frac{1}{3},1, \frac{1}{3}, 0 } \big)$ using 
\begin{equation}
 R_{\frac{1}{3},1, \frac{1}{3}, 0 }(z)= \frac{9 z^{3}}{2} - \frac{55 z}{8} - \frac{9}{8}.
\end{equation}
\begin{figure}[htb]
\centering
\includegraphics[width=0.85\textwidth]{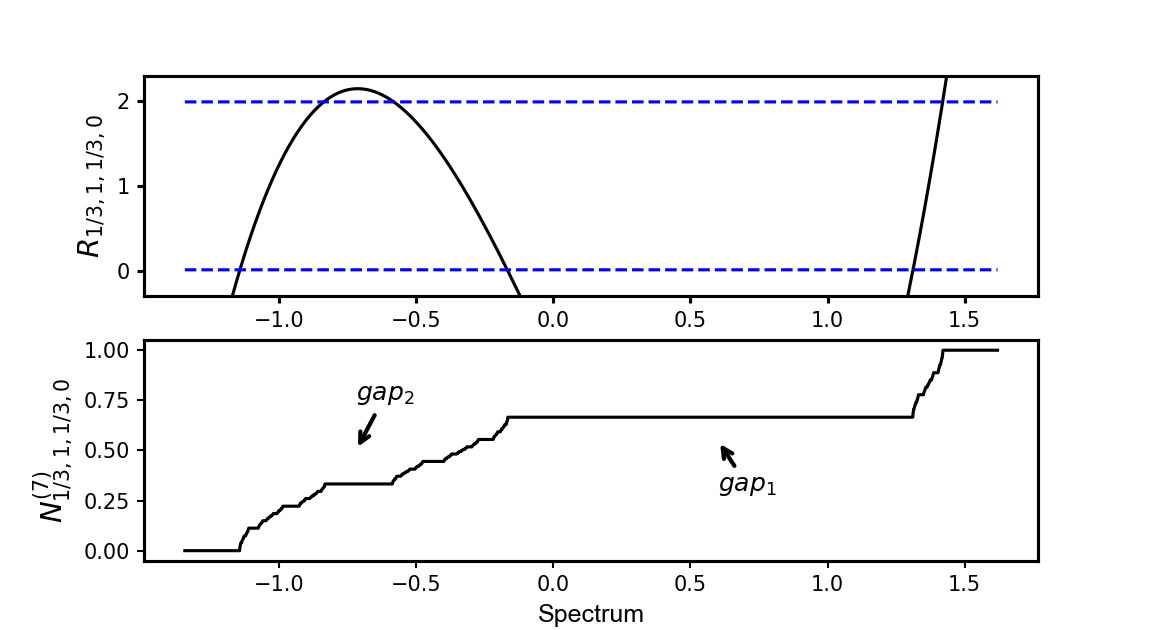}
\vspace{0cm}
\caption{(Top) The spectral decimation function $R_{\frac{1}{3},1, \frac{1}{3}, 0 }$ is plotted. The dashed  lines represent the cutoffs at $y=0$ and $y=2$. (Bottom)  The  integrated density of states $N^{(l)}_{\frac{1}{3},1, \frac{1}{3}, 0 }$  is plotted for level $l=7$. The dashed cutoff lines and spectral decimation function are used to locate the spectral gaps, which coincide with the indicated plateaus of the integrated density of states.}
\label{fig:GeneratinggapsAMOfirstTwoGaps}
\end{figure}
To locate the first two spectral gaps of $ \hamilton_{\frac{1}{3},1, \frac{1}{3}, 0 } $, we note that $\sigma(\Delta_{\frac{1}{3}}) \subset [0,2]$. By Theorem \ref{thm:SpectralSimiAMOandpqModelnew}, we obtain
\begin{align*}
z \in \sigma \Big( \hamilton_{\frac{1}{3},1, \frac{1}{3}, 0 }  \Big) \quad \Rightarrow \quad   R_{\frac{1}{3},1, \frac{1}{3}, 0 }(z) \in [0,2],
\end{align*}
or equivalently 
\begin{align*}
R_{\frac{1}{3},1, \frac{1}{3}, 0 }(z) \notin [0,2] \quad \Rightarrow \quad z \in \rho \Big( \hamilton_{\frac{1}{3},1, \frac{1}{3}, 0 }  \Big).
\end{align*}
Plotting the spectral decimation function $R_{\frac{1}{3},1, \frac{1}{3}, 0 }$ with both cutoffs $y=0$ and $y=2$ in Figure \ref{fig:GeneratinggapsAMOfirstTwoGaps}, generates the first two spectral gaps $gap_1$ and $gap_2$. By Proposition \ref{prop:PropertiesOfR}, we know that $z \in \sigma \big(\hamilton_{\frac{1}{3},1, \frac{1}{3}, 0 }^{(1)}\big) \bigcup \sigma \big(\hamilton_{\frac{1}{3},1, \frac{1}{3}, 0 }^{(1),D} \big)$ if and only if $R_{\frac{1}{3},1, \frac{1}{3},0}(z) \in  \{0,2\}$.
The eigenvalues of $\hamilton_{\frac{1}{3},1, \frac{1}{3}, 0 }^{(1)}$ are listed in Table  \hyperref[tab:table1]{1} and we denote the eigenvalues $\sigma \big(\hamilton_{\frac{1}{3},1, \frac{1}{3}, 0 }^{(1),D} \big) = \{-\frac{5}{6}, -\frac{1}{6} \}$  by $\lambda^{(1),D}_1 =-\frac{5}{6}$ and  $\lambda^{(1),D}_2 = -\frac{1}{6}$. This gives
\begin{align}
\label{eq:orderingEigenvalues}
\lambda^{(1)}_1 \leq  \lambda^{(1),D}_1   \leq \lambda^{(1)}_2 \leq    \lambda^{(1),D}_2  \leq \lambda^{(1)}_3 \leq \lambda^{(1)}_4.
\end{align}
The spectrum of $\hamilton_{\frac{1}{3},1, \frac{1}{3}, 0 }$ is then contained in the complement (in $\rr$) of the following set
\begin{align*}
(-\infty, \lambda^{(1)}_1) \cup (\lambda^{(1),D}_1,   \lambda^{(1)}_2) \cup (\lambda^{(1),D}_2, \lambda^{(1)}_3) \cup (\lambda^{(1)}_4,\infty),
\end{align*}
where $gap_1=(\lambda^{(1),D}_1,   \lambda^{(1)}_2) $ and $gap_2 = (\lambda^{(1),D}_2, \lambda^{(1)}_3)$. 
\begin{figure}[htb]
\centering
\includegraphics[width=0.9\textwidth]{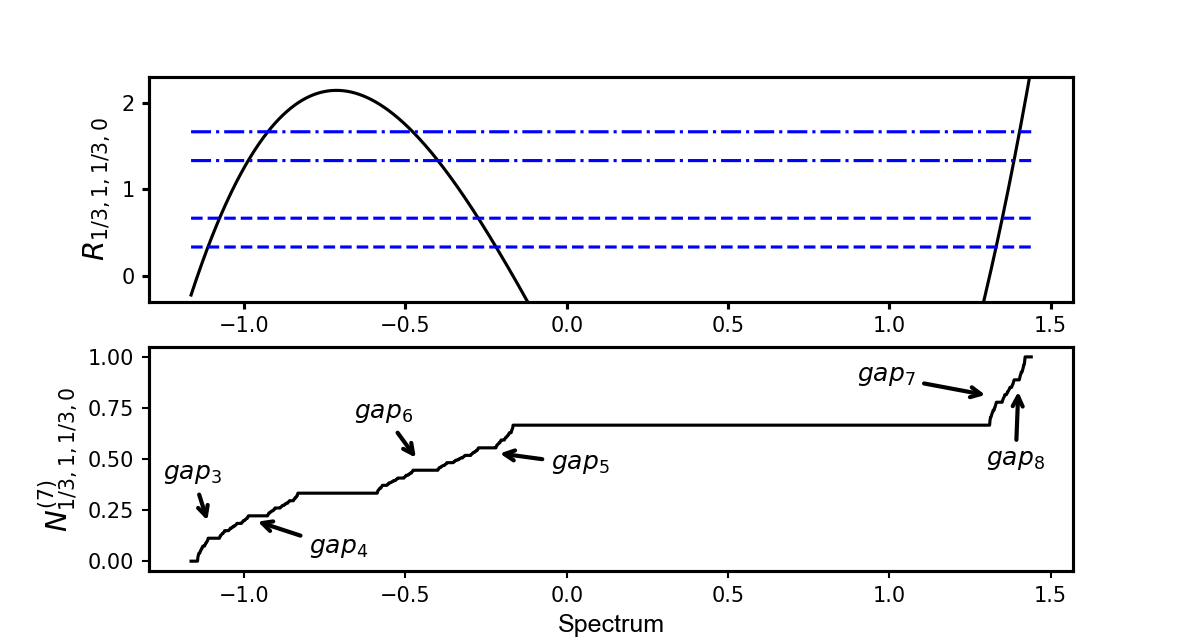}
\vspace{0cm}
\caption{(Top) The spectral decimation function $R_{\frac{1}{3},1, \frac{1}{3}, 0 }$ is plotted. The dashed  lines represent the cutoffs at $y=\frac{1}{3}$ and $y=\frac{2}{3}$ and the dash-dot lines represent the cutoffs at $y=\frac{4}{3}$ and $y=\frac{5}{3}$.  (Bottom)  The  integrated density of states $N^{(l)}_{\frac{1}{3},1, \frac{1}{3}, 0 }$  is plotted for level $l=7$. The cutoff lines and spectral decimation function are used to locate the spectral gaps, which coincide with the indicated plateaus of the integrated density of states.}
\label{fig:GeneratinggapsAMOnext6Gaps}
\end{figure}
To generate the next spectral gaps, we proceed similarly and note that
\begin{equation}
\sigma(\Delta_{\frac{1}{3}}) \subset \big[0,\frac{1}{3}\big] \cup \big[\frac{2}{3}, \frac{4}{3}\big] \cup \big[\frac{5}{3}, 2\big] 
\end{equation}
where $\sigma(\Delta^{(1)}_{p})=\{0,\frac{1}{3},\frac{5}{3},2\}$ and 
$\sigma(\Delta^{(1),D}_{p})=\{\frac{2}{3},\frac{4}{3}\}$ (with Dirichlet boundary conditions). Hence,
\begin{align*}
R_{\frac{1}{3},1, \frac{1}{3}, 0 }(z) \in \Big( \frac{1}{3},\frac{2}{3} \Big) \cup \Big( \frac{4}{3},\frac{5}{3} \Big) \quad \Rightarrow \quad z \in \rho \Big( \hamilton_{\frac{1}{3},1, \frac{1}{3}, 0 }  \Big).
\end{align*}
Plotting the spectral decimation function $R_{\frac{1}{3},1, \frac{1}{3}, 0 }$ with both cutoffs $y=\frac{1}{3}$, $\frac{2}{3}$ and $y=\frac{4}{3}$, $y=\frac{5}{3}$ in Figure \ref{fig:GeneratinggapsAMOnext6Gaps}, generates the next six spectral gaps.

\subsection{Gap labeling}\label{subsec:numerical}

\begin{figure}[!htb]
  \begin{minipage}[b]{0.5\linewidth}
    \centering
    \includegraphics[width=1.1\linewidth]{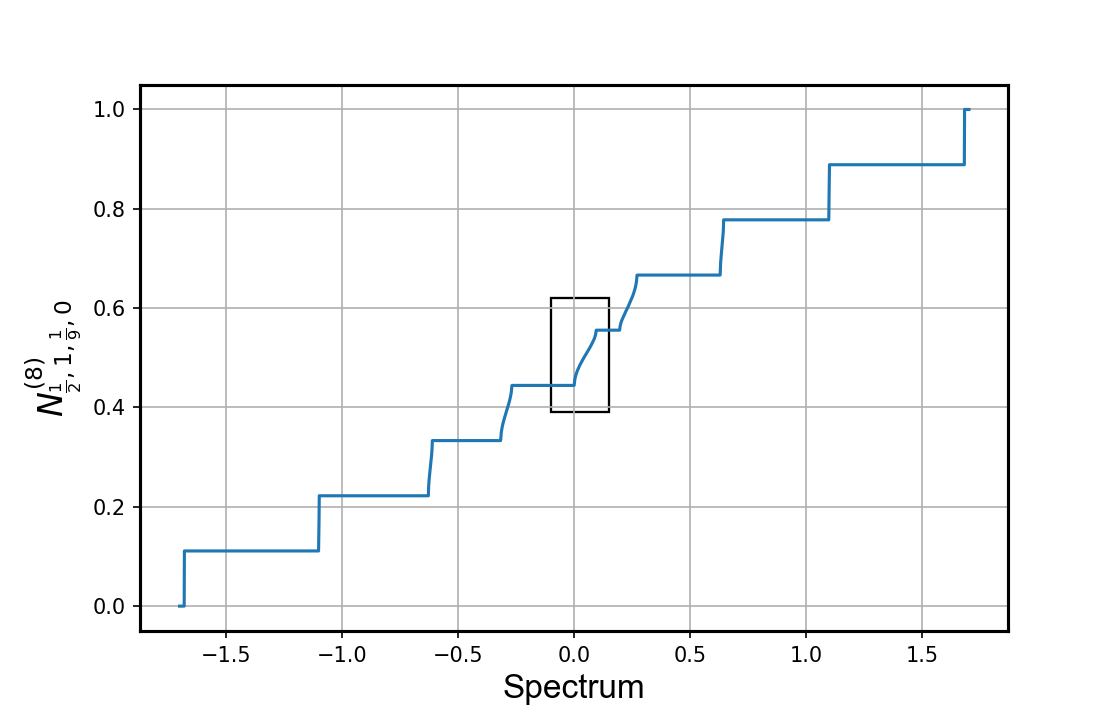} 
  \end{minipage}
  \begin{minipage}[b]{0.5\linewidth}
    \centering
    \includegraphics[width=1.1\linewidth]{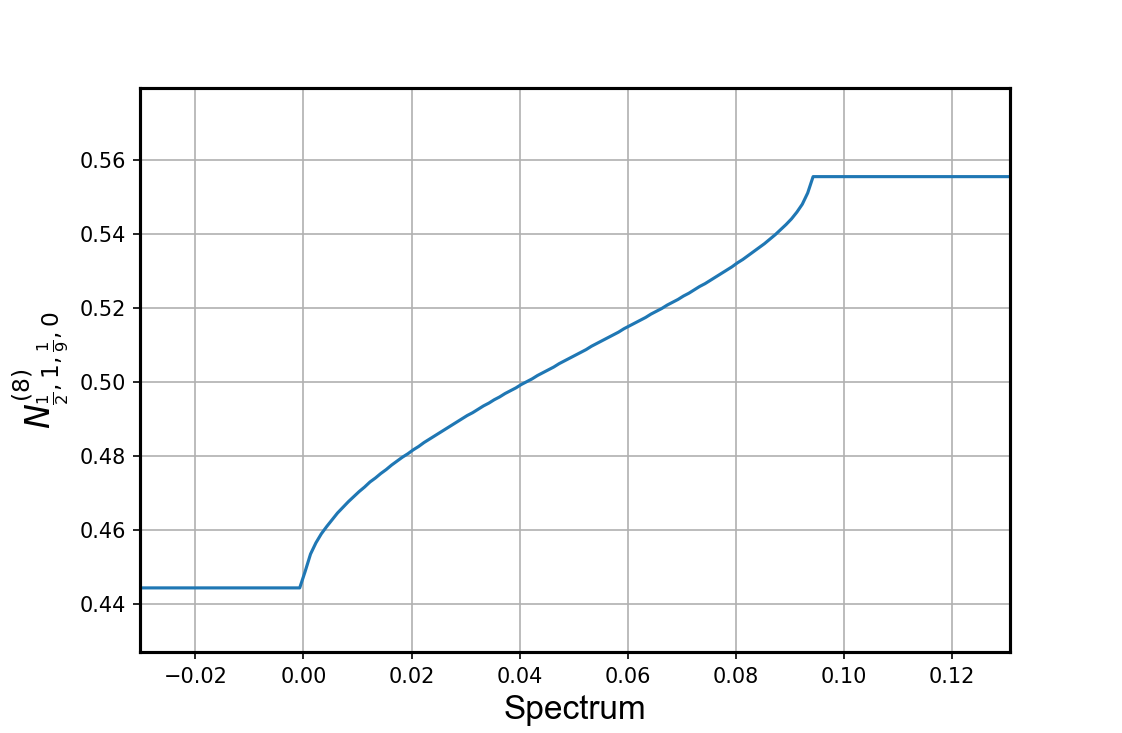} 
  \end{minipage} 
  \caption{ (Left) Numerical computation of the integrated density of states for  $\hamilton_{p,\beta, \alpha, \theta }^{(l)}$. The computations are done for level $l=8, p=\frac{1}{2}, \beta=1,\alpha=\frac{1}{9},\theta=0$. (Right) A resized version of the box in the left-hand side figure is displayed. It shows a copy of Figure   \ref{fig:IDSforLaps} (left) rescaled by $\frac{1}{9}$.}
  \label{fig:IDSAMOalpha1over9p1over2}
\end{figure}

\begin{figure}[!htb]
  \begin{minipage}[b]{0.5\linewidth}
    \centering
    \includegraphics[width=1.1\linewidth]{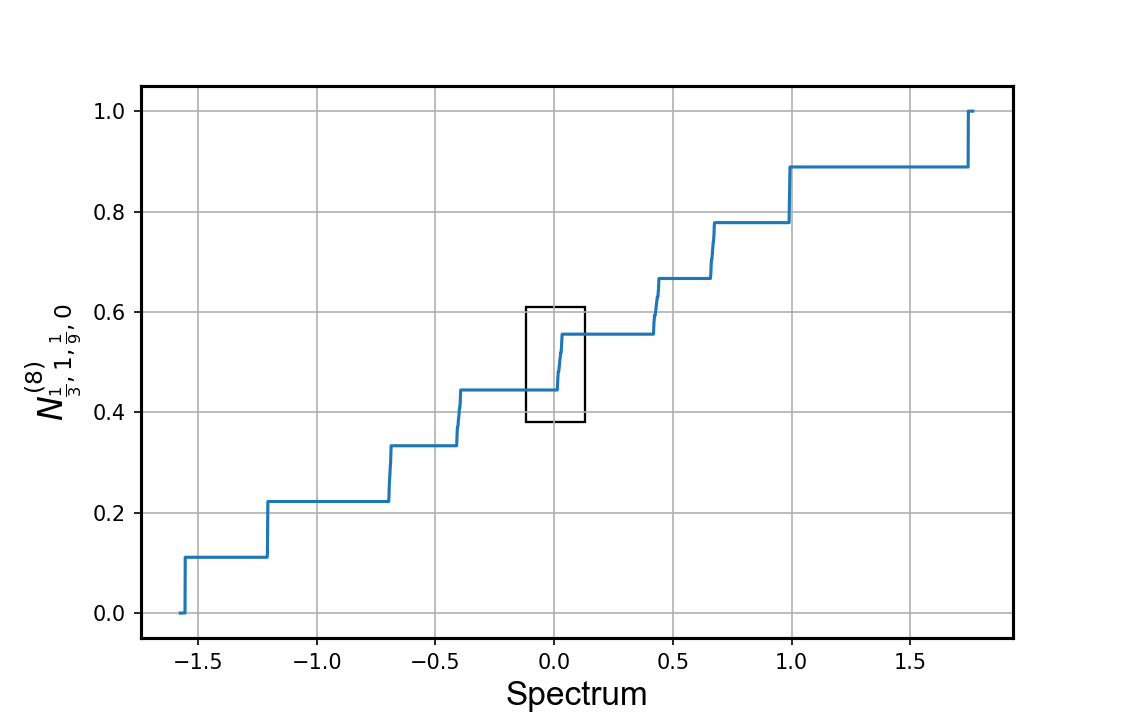} 
  \end{minipage}
  \begin{minipage}[b]{0.5\linewidth}
    \centering
    \includegraphics[width=1.1\linewidth]{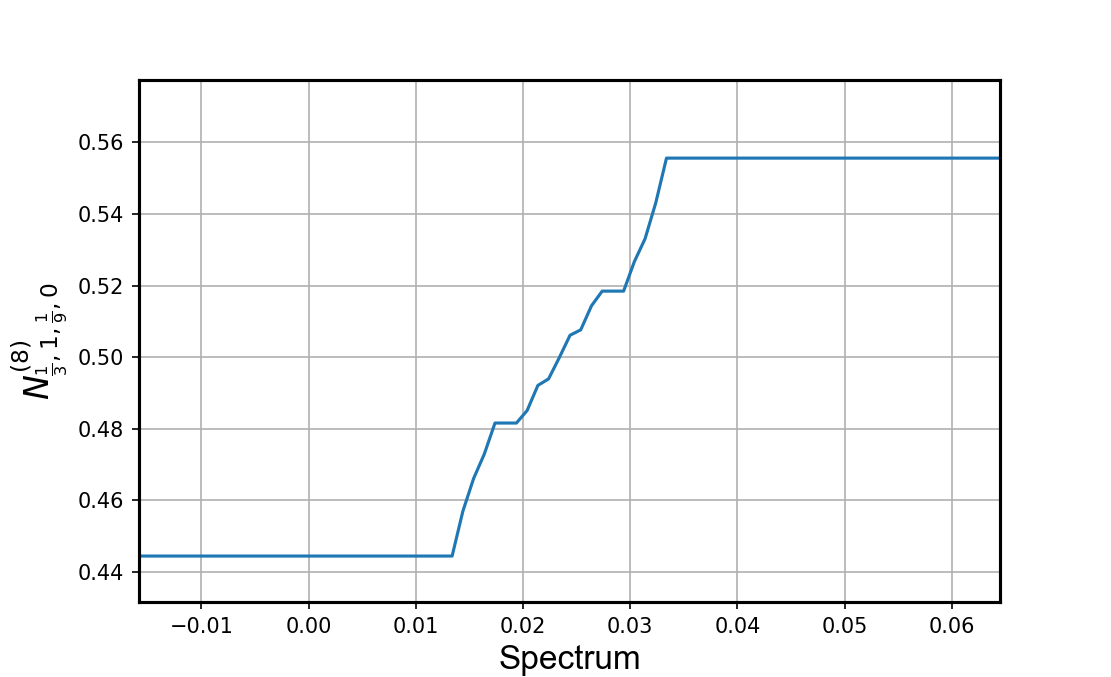} 
  \end{minipage} 
  \caption{ (Left) Numerical computation of the integrated density of states for  $\hamilton_{p,\beta, \alpha, \theta }^{(l)}$. The computations are done for level $l=8, p=\frac{1}{3}, \beta=1,\alpha=\frac{1}{9},\theta=0$. (Right) A resized version of the box in the left-hand side figure is displayed. It shows a copy of Figure   \ref{fig:IDSforLaps} (right) rescaled $\frac{1}{9}$.}
  \label{fig:IDSAMOalpha1over9p1over3}
\end{figure}

Figures \ref{fig:IDSAMOalpha1over9p1over2} and  \ref{fig:IDSAMOalpha1over9p1over3} give a numerical  illustration of~\eqref{eq:almostSelfSimiId}). We recall that  the spectral decimation function $R_{\frac{1}{2},1, \frac{1}{9},0 }$ is a polynomial of degree $9$. As such,  Figure \ref{fig:IDSAMOalpha1over9p1over2} shows that on each range of the nine branches $S_i(\Delta_{\frac{1}{2}})$, we have a copy of Figure~\ref{fig:IDSforLaps} (left) rescaled by $\frac{1}{9}$, according to Corollary \ref{cor:scalingCopies}. This case corresponds to a periodic Jacobi matrix, and the spectrum consists of nine spectral bands. As expected from Corollary~\ref{coro:GapLabel}, the set of spectral gap labels is
\begin{align}
\mathscr{G} \mathscr{L}(\hamilton_{\frac{1}{2},1, \frac{1}{9}, 0 }) = \Big\{ 0, \frac{1}{9}, \frac{2}{9}, \dots , 1   \Big\}.
\end{align}

Similarly, the spectral decimation function $R_{\frac{1}{3},1, \frac{1}{9},0 }$ is a polynomial of degree $9$. As such, Figure~\ref{fig:IDSAMOalpha1over9p1over3} shows that on each range of the nine branches $S_i(\Delta_{\frac{1}{3}})$, we have a copy of Figure~\ref{fig:IDSforLaps} (right) rescaled by $\frac{1}{9}$. In particular, this highlight the Cantor set structure inherited from  $\sigma(\Delta_{\frac{1}{3}})$ and the set of spectral gap labels is deduced from Corollary \ref{coro:GapLabel}.

\FloatBarrier

\section{Connections approaches of    B\'{e}llissard and Bessis-Geronimo-Moussa}\label{sec:Bellissard}
 
{   Our work is partially motivated by Bellissard's studies on Hamiltonians describing the motion of a particle in quasicrystals. More specifically, his construction of  a large class of Hamiltonians with Cantor spectra starting from Jacobi matrices and their  associated Julia sets, see  \cite{BellissardRenormalizationGroup1992,Bellissard85}. To draw the parallel with our work we make the following observations. Given a polynomial  $P(z)$, let $\mathcal{J}(P)$ be the corresponding (compact) Julia set, which, under some assumptions on $P(z)$ is a completely disconnected set \cite{Fatou1919,Julia1918}. In addition, let  $\mu$  be the balanced invariant measure of $\mathcal{J}(P)$ and consider the Hilbert space $\mathcal{H}=L^2(\mathcal{J}(R), d\mu)$. The multiplication operator $H$ associated with the identity function $f(x)=x$ in $\mathcal{H}$ is bounded, self adjoint, has the Cantor set $\mathcal{J}(P)$  as spectrum. Furthermore, $\mu$ is  the spectral measure of $H$ leading to a singular spectrum. Because the linearly independent set $\{x \in \mathcal{J}(P) \to x^n \in \complex\}$ generates a dense linear subspace of  $\mathcal{H}$, the operator $H$ can be represented  as a semi-infinite Jacobi matrix\cite{Barnsley1983,Barnsley1985}. Moreover, $H$ satisfies a \textit{renormalization group equation}
 \begin{align}
 \label{eq:RenoGroBelli}
 D \big( zI-H \big)^{-1} D^{\ast} =\frac{P'(z)}{deg(P)} \big( P(z)I-H \big)^{-1}
\end{align}
where the partial isometry $D$ and its adjoint $D^{\ast}$ are given in~\cite[Theorem 1]{BellissardRenormalizationGroup1992}.

The connection between our work and Bellissard's original ideas as elaborated above begins by using defining a probabilistic Laplacian (the so-called \textit{pq-model}) on the integers half-line $ \mathbb{Z}_+$, regarded as a hierarchical or substitution graph with $G_1$ in Figure~\ref{fig:G1_initialNeumann} (right) as its basic building block.
The graph $G_1$ determines the spectral decimation function $R_{\Delta_{p}}(z)$ in~\eqref{eq:specDeciFuncLap}, a polynomial which plays the role of the polynomial $P(z)$  appearing in Bellissard's approach. On the one hand,  each Laplacian $\Delta_p$ represents an example of a semi-infinite  Jacobi matrix that, in similarly to the operators in Bellissard's construction, has a spectrum that coincides with $\mathcal{J}(R_{\Delta_{p}})$: the Julia set of the polynomial $R_{\Delta_{p}}(z)$. On the other hand, the balance measure of $\mathcal{J}(R_{\Delta_{p}})$ plays the role of the spectral measure in Bellissard's approach and is also the density of states in our context, see Proposition~\ref{prop:IDSofpLap}. In a forthcoming work \cite{BaluMograbyOkoudjouTeplyaevJacobi2021}, we generalize the substitution rule in Definition~\ref{def:finiteGraphApproxNeumann} leading to a multiple-parameter families of probabilistic Laplacians whose spectral properties are investigated with similar tools as the aforementioned  \textit{pq-model}.

We note that the spectra of the self-similar almost Mathieu operators we introduced in this paper are not necessarily given as the Julia sets of some polynomials. Instead,  as proved in Theorem~\ref{thm:SpectralSimiAMOandpqModelnew} these  spectra are preimages of Julia sets under certain polynomials. In \cite{BaluMograbyOkoudjouTeplyaevJacobi2021}, these results are extended to a class of Jacobi operators, for which we establish  a renormalization group equation, see \cite[Theorem 4.6.]{BaluMograbyOkoudjouTeplyaevJacobi2021}. For example, when reduced to the \textit{pq-model}, this renormalization group equation for resolvent, see \cite{Luke10,Luke12,TeplyaevInfiniteSG1998,malozemovteplyaev2003}, takes the form 
\begin{align}
\label{eq:RenoGrEqOURapproach}
U^{\ast}\big(z I- \Delta_{p} \big)^{-1}U=\frac{(z-1)^2-p^2}{p(1-p)}\Big( R_{\Delta_{p}}(z) I - \Delta_{p} \Big)^{-1},
\end{align}
where $R_{\Delta_{p}}(z)$ is given in~\eqref{eq:specDeciFuncLap}, and $U$ and $U^{\ast}$ are defined in  \cite{ChenTeplyPQmodel2016}. For comparison with~\eqref{eq:RenoGroBelli}, we compute 
\begin{align}
\frac{\frac{d}{dz}R_{\Delta_{p}}(z) }{deg(R_{\Delta_{p}})} = \frac{(z-1)^2-\big(\frac{1-p(1-p)}{3}\big)}{p(1-p)},
\end{align}
which coincides with the factor on the right-hand side of equation (\ref{eq:RenoGrEqOURapproach}) if and only if $p=\frac{1}{2}$.
Thus, our results is related to \cite[Theorem 2.2]{BGM-1988}, although we do not rely on \cite{BGM-1988} as we consider a model which allows us to produce more explicit computation of the spectrum for operators with potential.

}

\section{Conclusions}

In this paper we  introduce and study a fractal version of the almost Mathieu operators $\hamilton_{p,\beta, \alpha, \theta }$. We propose a new adaptation of the  method of spectral similarity to analyze their spectral properties. Our main conclusions are the following. 
\begin{enumerate}
\item  Theorem \ref{thm:SpectralSimiAMOandpqModelnew} presents a useful algebraic tool to relate the spectrum of the almost Mathieu operators  $\hamilton_{p,\beta, \alpha, \theta }$ to that of a family of self-similar Laplacians $\Delta_p$. Our results are established when the parameter $\alpha$ belongs to the dense set $\{\tfrac{k}{3^n}, \, k=0,1, 2, \hdots, 3^n-1\}_{n=1}^l$ where $l\geq 1$. Note that in the classical case, corresponding to $p=1/2$ in our formulation, many important results   are also obtained  for $\alpha$ irrational, but in the fractal setting the methods for irrational $\alpha$ are not developed yet. 

 The method of spectral similarity is applicable in many situations, and in a forthcoming work \cite{BaluMograbyOkoudjouTeplyaevJacobi2021} we develop a general framework by working with  a large class of Jacobi operators. We call these operators piecewise centrosymmetric Jacobi operators  \cite{centrosymmetric1,centrosymmetric2,centrosymmetric3}. 
 In this general setting   the spectral decimation function can be computed using the theory of   orthogonal polynomials associated with the aforementioned Jacobi matrix. As a result, the spectral decimation function in the generalizations of Theorem \ref{thm:firstResultFiniteGraphs} can be computed using a three-term recurrence relation, which provides a simple procedure to show that the spectral decimation function is a polynomial of a specific degree with properties that can be controlled.

\item Our methods allow to compute explicitly the density of states. In particular, we proved in theorem \ref{thm:IDSofAMO} an explicit formula connecting   the density of states of  $\hamilton_{p,\beta, \frac{k}{3^n}, 0 }$ by identifying it with the weighted preimages of the balanced invariant measure on the Julia set of the polynomial  $\mathcal{J}(R_{\Delta_{p}})$. This approach can be generalized to many other situations, and can be verified numerically, see  Section~\ref{sec:examplesappl}.  

\item In our particular situation  we are able to conclude  that the operators  $\hamilton_{p,\beta, \frac{k}{3^n}, 0 }$ have  singularly  continuous spectrum when $p\neq \frac{1}{2}$ because of the previous recent work \cite{ChenTeplyPQmodel2016}  that the spectrum is singularly continuous for  $\Delta_p$. This result requires a detailed analysis of a certain dynamical system describing the behavior of the generalized eigenfunction.

\item A particular novelty of our results is that we develop the spectral analysis of a self-similar Laplacian with a quasi-periodic potential. In our work we have made the essential steps towards the Fourier analysis on one-dimensional self-similar structures   following the general approach developed by Strichartz et al  \cite{Bob1989HarmonicAnalysis,Strichartz2003FractafoldsBO}.  In certain particular situations this allows to consider classical and quantum wave prorogation on fractal and other irregular structures \cite{Akkermans2013StatMechFractals,Akkermans2014WavePropFractal,wave17}.

\item Using these Harmonic Analysis tools, our work introduces a direct approach to the  gap labeling for a self-similar Laplacian with a potential. In particular, this   direct approach   for gap labeling is complementing \cite{be1,be2,be3}. In our case we do not make use of the   dual of a group acting on the fractal lattice, and observe that gap labels of the form $\frac{j}{3^n}$ are consistent  with the self-similar quasi-periodic structure where renormalization acts by the dilation of the space  by $3^n$. This should be contrasted with the gap labeling for fractal structures with more complicated topological structure currently under investigation in \cite{BubbleDiamond2021}, {including the classical Sierpinski gasket and a new model of the bubble diamond fractals}. In general, on a certain class of self-similar structures the gaps  are labeled by the values of the 
 integrated density of states of the Laplacian with values $\frac{j}{C^n}$ where $C$ is the topological degree of the self-covering map of the fractal limit space. 
 Thus, our work  sets the stage for considering the Bloch theorem, noncommutative Chern characters and fractal-based quantum Hall systems, see   \cite{Marcolli2006}, as well as \cite{Benameur-Mathai}. On fractal spaces this is an open problem that has not been previously considered in the literature besides the recent work \cite{Akkermans2021}. 
 
 \item Our work is connected to several lines of investigations in mathematics and physics, including the topics highlighted at the recent workshop 
 \href{https://alexander-teplyaev.uconn.edu/quasi-periodic-spectral-and-topological-analysis/}{Quasi-periodic spectral and topological analysis}, 
 and in particular the work of 
 E.~Akkermans \cite{akke17,akkermans2007mesoscopic,ovdat2021breaking,Akkermans2014WavePropFractal,Akkermans2013StatMechFractals,Akkermans2021},
 D.~Damanik   \cite{Damanik2021,Damanik2019,Damanik2015,Damanik2004,Damanik1999},
 S.~Jitomirskaya \cite{jitomirskaya_metal-insulator_1999,marx_jitomirskaya_2017,jito19,jito20,jito21,avila_ten_2009}, and 
 E.~Prodan \cite{Prodan2016a,Prodan2016b,Prodan2016c,rosa2021topological,ni2019observation} et al.

\end{enumerate}

%

\subsection*{Acknowledgments} The work of G. ~Mograby was supported by ARO grant W911NF1910366.  K.~A.~Okoudjou was partially supported by ARO grant W911NF1910366 and the National Science Foundation under Grant No. DMS-1814253. A.~Teplyaev was partially supported by NSF DMS grant 1613025 and by the Simons Foundation.

\bibliographystyle{plain}
\bibliography{NewAMO-Ref}

\end{document}